\documentclass[12pt]{article} 
\usepackage{abstract}
\usepackage[utf8]{inputenc} 
\usepackage{amsxtra} 
\usepackage{amsthm}
\usepackage{amssymb}
\usepackage{amsfonts}
\usepackage{graphicx}
\usepackage{rotating}
\usepackage{color}
\usepackage{relsize}
\usepackage{epsfig}
\usepackage{subcaption} 
\usepackage{verbatim}
\usepackage[printonlyused]{acronym} 
\usepackage{setspace}
\usepackage{epstopdf}
\usepackage{authblk}
\usepackage{sectsty}
\usepackage[colorinlistoftodos]{todonotes}
\usepackage[colorlinks=true,allcolors=black]{hyperref}
\usepackage{textcomp,marvosym}
\usepackage{soul}
\usepackage{algorithm}
\usepackage{algorithmicx}
\usepackage{algpseudocode}
\usepackage{url}
\usepackage{hyperref}
\usetikzlibrary{shapes,arrows,automata,positioning}
\usepackage{mathtools}

\sectionfont{\fontsize{12}{15}\selectfont} 
\subsectionfont{\fontsize{10}{12}\selectfont}

\usepackage{geometry} 
\definecolor{blue}{RGB}{0, 0, 150}
\definecolor{green}{RGB}{0,150,0}
\definecolor{red}{RGB}{200, 0, 0}
\definecolor{black}{RGB}{0, 0, 0}
     
\renewcommand{\l}{\left}
\renewcommand{\r}{\right}

\newcommand{\red}[1]{{\color{red}{#1}}}

\newcommand{\dd}{\text{\,d}}
\definecolor{Purpleee}{RGB}{91,40,224}

\DeclareMathOperator*{\argmin}{arg\,min}

\begin{document}

\title{Regression-based projection for learning Mori--Zwanzig operators}
\author[1,$\dagger$,*]{Yen Ting Lin}
\author[1,$\dagger$]{Yifeng Tian}
\author[3]{Danny Perez}
\author[2]{Daniel Livescu}

\affil[1]{Information Sciences Group, Computer, Computational and Statistical Sciences Division (CCS-3), Los Alamos National Laboratory, Los Alamos, NM 87545, USA}	
\affil[2]{Computational Physics and Methods Group, Computer, Computational and Statistical Sciences Division (CCS-2), Los Alamos National Laboratory, Los Alamos, NM 87545, USA}
\affil[3]{Physics and Chemistry of Materials Group, Theoretical Division (T-1), Los Alamos National Laboratory, Los Alamos, NM 87545, USA}

\affil[$\dagger$]{These authors contributed equally.}

\affil[*]{Corresponding author: yentingl@lanl.gov}
\date{\today}

	\maketitle
    \begin{abstract}
{  We propose to adopt statistical regression as the projection operator to enable data-driven learning of the operators in the Mori--Zwanzig formalism. We present a principled method to extract the Markov and memory operators for any regression models. We show that the choice of linear regression results in a recently proposed data-driven learning algorithm based on Mori's projection operator, which is a higher-order approximate Koopman learning method. We show that more expressive nonlinear regression models naturally fill in the gap between the highly idealized and computationally efficient Mori's projection operator and the most optimal yet computationally infeasible Zwanzig's projection operator. We performed numerical experiments and extracted the operators for an array of regression-based projections, including linear, polynomial, spline, and neural-network-based regressions, showing a progressive improvement as the complexity of the regression model increased. Our proposition provides a general framework to extract memory-dependent corrections and can be readily applied to an array of data-driven learning methods for stationary dynamical systems in the literature.   \vspace{6pt}\\
    {\bf Keywords:} Mori--Zwanzig formalism, Koopman representation, nonlinear projection operators, data-driven learning, regression, neural networks \vspace{6pt}\\
    {\bf AMS subject classifications:}  37M19, 37M99, 46N55, 65P99, 82C31}
	\end{abstract}
\vspace{10pt}

\section{Introduction}
More than half a century ago, Mori \cite{mori1965transport} and Zwanzig \cite{zwanzig1973nonlinear} developed a mathematically exact formalism for constructing reduced-order models for dynamical systems in non-equilibrium statistical mechanics. In this context, the typical degrees of the freedom of a system are of the order of Avogadros number. In contrast to tracking the dynamics of each degree of freedom, reduced-order models (also referred to as the coarse-grained models in the context of statistical physics) are developed to describe the dynamics of a relatively small number of variables  of  interest. In the jargon of model coarse-graining, this set of dynamic variables is referred to as the ``resolved'' variables, and the rest degrees of freedom in the system is referred to as the ``unresolved'' ones. Reduced-order models are particularly useful in describing the emergent phenomena at the mesoscopic or macroscopic scales. They are also more computationally efficient to simulate, and thus play an essential role in bridging multiscale models. 

The key difficulty in constructing reduced-order models is the \emph{closure problem}: one must quantify the effect of the unresolved variables before self-contained evolutionary equations of the resolved variables can be prescribed. To solve the closure problem, Mori and Zwanzig adopted an elegant approach which uses functional projections. Formally, a pre-specified projection operator is used to map any function which depends on both resolved and unresolved variables to a function that only depends on the resolved variables. As the latter does not depend on the unresolved information, one can then prescribe a set of evolutionary equations for the resolved variables. The dynamical system in terms of only the projected functions is of course only an approximation of the true dynamics. The Mori--Zwanzig (MZ) formalism further quantifies the error of the approximation and derives how the error would propagate into the future by the dynamics. The central result of the Mori--Zwanzig formalism is the Generalized Langevin Equation, which stipulates that the evolutionary equations of the resolved variables do not only depend on their instantaneous values but also their past history. As a mathematically exact formula, the memory dependence, which depends critically on the choice of the resolved variables and the projection operator, quantifies the interactions between the resolved and unresolved degrees of freedom.

Theoretically speaking, there are infinitely many choices of the projection operators. In the literature, the commonly chosen projections include (1) Mori's linear projection \cite{mori1965transport}, (2) finite-rank projection operator \cite{chorin2002OptimalPredictionMemory}, which is Mori's projection with a set of orthonormal basis functions, (3) projection by assigning unresolved variables to zero \cite{parishDynamicSubgridScale2017,parishNonMarkovianClosureModels2017,stinis19,stinis21}, and (4) Zwanzig's fully nonlinear projection \cite{zwanzig1973nonlinear,chorin2002OptimalPredictionMemory}. It is generally difficult to derive analytic expressions for the Mori--Zwanzig formalism. Several recent studies \cite{okamura2006validity, mori2007dynamic, meyer2020non, maeyama2020extracting, meyer2021numerical,lin2021DataDrivenLearningMori,tian2021DatadrivenLearningMori} have established that with Mori's linear projection operator, it is possible to adopt a data-driven approach to learn the Mori--Zwanzig operators using the time series of the resolved dynamics. In addition, it was shown \cite{lin2021DataDrivenLearningMori} that these methods provide higher-order and memory-dependent corrections to existing data-driven learning of the approximate Koopman operators \cite{rowley2009SpectralAnalysisNonlinear,schmid2010DynamicModeDecomposition,Schmid2011,williams2015DataDrivenApproximation}.

In spite of being a consistent theoretical framework to the approximate Koopman learning, Mori's projection operator is greatly limited by its linear nature. That is, after the projection, all functions must be expressed as linear combinations of a set of \emph{a priori} specified resolved variables. The quality of the prediction thus solely depends on the choice of the resolved variables. To the authors' best knowledge, there is no principled way to select an optimal set of resolved variables for general dynamical systems. With a non-optimal choice, the difference between the projected dynamics and the true dynamics can be unsatisfactorily large \cite{lin2021DataDrivenLearningMori}. The same problem also manifests itself in the approximate Koopman learning framework \cite{williams2015DataDrivenApproximation}. Furthermore, Mori's projection operator is not compatible with many nonlinear closure schemes \cite{Durbin18,Majda06}. Neither can we generalize the Mori's restrictive linear projection operator to those data-driven learning methods based on modern machine learning (ML) \cite{Brunton16SINDy,liExtendedDynamicMode2017,Yeung2017LearningDN,luschDeepLearningUniversal2018,wehmeyerTimelaggedAutoencodersDeep2018}, with which one relies on the neural networks (NNs) to identify the nonlinear relationship between the input (what we know at the present) and the output (what we want to predict in the future).

At the other end of the spectrum is the Zwanzig's construct \cite{zwanzig1973nonlinear,chorin2002OptimalPredictionMemory} which uses conditional expectation as a projection operator. Although the conditional expectation is the optimal choice of the projection operator \cite{chorin2002OptimalPredictionMemory}, in practice, it is challenging to estimate the conditional expectations with a finite set of data, especially for high-dimensional systems that are not fully-resolved. 

To address the aforementioned difficulties, we propose two novel ideas that enable data-driven learning for the Mori--Zwanzig operators. First, we propose a concept that \textit{any regression analysis} can be treated as \textit{the} projection operator in the Mori-Zwanzig formalism. Secondly, we identify that the Generalized Fluctuation-Dissipation (GFD) relationship, a self-consistent condition that relates the Mori--Zwanzig memory kernels and the error of the memory-dependent dynamics, \textit{should be} used to learn the memory kernels recursively. Building upon this concept, we have developed a systematic approach for learning Mori-Zwanzig memory kernels, using snapshots of resolved variables. The resulting procedure is surprisingly simple, principled, and intuitive, despite a complicated operator-algebraic derivation. It is important to note that the memory kernels depend on the choice of the projection operator in the Mori-Zwanzig formalism. As such, our procedure delivers memory kernels that are specific to the choice of regression model. We will demonstrate that with a linear regression model, the learned memory kernels converge nicely to the outcome of our recently proposed data-driven learning for Mori's projector \cite{lin2021DataDrivenLearningMori}. Moreover, we will elaborate on how seemingly similar regression models can be associated with distinct projection operators, and how differentiating such a nuance could pave the way to enable learning the Mori--Zwanzig memory kernels for reduced-order models with nonlinear closure schemes. 

Before the technical presentation, it may be beneficial to provide a high-level characterization that contrasts our proposition to other statistical and machine-learning models motivated by Mori--Zwanzig's memory-dependent formulation. By doing so, we hope to provide a broader perspective on our work and the contributions it makes to the field. The key idea to differentiate our work to other models lies in the fact that Mori--Zwanzig's formalism is memory-dependent, but not all memory-dependent formulation is Mori--Zwanzig. The memory kernels in Mori--Zwanzig formalism are special because they \textit{must} satisfy the self-consistent GFD. As such, it would be a logical fallacy to call \textit{any} memory-dependent dynamics as a Mori--Zwanzig model without enforcing GFD in the model structure, or without checking the GFD after learning. Unfortunately, many existing models fell into this logical fallacy. For example, \cite{weinanE} directly jumped from the Mori--Zwanzig formalism to the Recurrent Neural Network with Long Short Term Memory (RNN-LSTM), claiming that the Mori--Zwanzig memory can be learned by black-box RNN-LSTM without enforcing or checking the GFD. The same gap exists in \cite{HARLIM2021109922} between the mathematical analysis based on Mori--Zwanzig formalism and the computational RNN-LSTM architecture. Although GFD was used as an argument to adopt finite-memory length truncation, it is not explicitly enforced or checked that the (trained) LSTM memory satisfies the GFD. With evidence presented \cite{weinanE} and \cite{HARLIM2021109922}, we cannot rule out that the RNN-LSTM learns other memory-dependent dynamics, for example, time delay embedding with which one augments the current state by the past history. Relevant data-driven delay embedded models include \cite{GILANI2021132829} which uses kernel methods and \cite{linDatadrivenModelReduction2021} which uses Wiener projection, which uses infinitely long past history. Both delay-embedding studies aim to \textit{bypass} quantification of the Mori--Zwanzig memory kernel, as they argued that optimal delay embedding would result in vanishing Mori--Zwanzig memories. A memory-dependent time-series analysis method termed NARMAX (nonlinear autoregressive moving average model with exogenous inputs \cite{chorinLu15}) also argues that its memory-dependence formulation resembles the Mori--Zwanzig formalism without establishing the GFD. However, it is later established in \cite{linDatadrivenModelReduction2021} that the memory-dependent formulation NARMAX is an approximation of the Wiener projection, which has zero Mori--Zwanzig memory. In contrast to these proposed methods, which postulate a memory-dependence in the statistical learning model without enforcing the GFD, our proposed approach has GFD built in. In fact, GFD is so critical that it turned out to the key equation to establish our proposed recursive procedure. To our best knowledge, such a procedure does not exist in the literature. 

To support the claim that our proposition can be applied to a wide range of data-driven models parametrized by regression, we will present numerical experiments on various nonlinear regression models, including polynomial and ridge regression, fully-connected neural networks (FCNN), and convolutional neural networks (CNN). It is not our intention to show that these simple regression models outperform existing data-driven and regression-based models. Rather, we would like to use these simple examples to illustrate the mechanism of the proposed procedure, the \text{learned Mori-Zwanzig (with enforced GFD) memories}, and how including the memory contribution can improve the predictions. We reemphasize that our proposition is not limited to these simple regression models. With the identification of a regression analysis as the Mori--Zwanzig projection operator, our proposed procedure can be applied to a wide array of data-driven models parametrized by regression (e.g., \cite{Brunton16SINDy,liExtendedDynamicMode2017,Yeung2017LearningDN,luschDeepLearningUniversal2018,wehmeyerTimelaggedAutoencodersDeep2018,QIAN2020132401}) for quantifying their Mori--Zwanzig memory kernels and for improving their predictions. Consequently, our proposition has the potential to significantly broaden the practical scope of the Mori-Zwanzig formalism.

\section{Background and terminology} \label{sec:background}

\subsection{Full dynamical systems} Following the notation in \cite{lin2021DataDrivenLearningMori}, we consider an autonomous and deterministic dynamical system in $\mathbb{R}^D$ following a continuous-time evolutionary equation 
\begin{equation}
    \frac{\dd}{\dd t} \boldsymbol{\phi}(t) = \mathbf{R}\left(\boldsymbol{\phi}\left(t\right)\right), \label{eq:ODE}
\end{equation}
where $\boldsymbol{\phi}$ is the state in the phase space  $\mathbb{R}^D$ and $\mathbf{R}:\mathbb{R}^D \rightarrow \mathbb{R}^D$ is the vector field of the system. We will assume that the system has a unique $\boldsymbol{\phi}(t)$ $\forall t \ge 0$. We will denote the solution of \eqref{eq:ODE} with the initial condition $\boldsymbol{\phi}_0 \in \mathbb{R}^D$ by $\boldsymbol{\phi}(t;\boldsymbol{\phi}_0)$. We will use the standard notation and denote scalars by symbols in the normal font, and vectors by symbols in bold.

\subsection{Observables} In reduced-order modeling, we aim to prescribe the evolutionary equations for $M < D$ real-valued functions of the state variable $\boldsymbol{\phi}$, $g_i:\mathbb{R}^D\rightarrow \mathbb{R}$, $i=1\ldots M$. We exclusively consider real-valued functions, but the theory can be straightforwardly generalized to complex-valued functions. We refer to these functions as the observables. In other context, such as material science or statistical mechanics, observables can also be referred to as descriptors, coarse-grained variables, or dynamic variables. We refer to an ``observation'' at time $t$ as applying the set of observables to the system's state $\boldsymbol{\phi}(t;\boldsymbol{\phi}_0)$, that is, a real-valued array $\l\{g_i \l(\boldsymbol{\phi}\l(t; \boldsymbol{\phi}_0\r)\r)\r\}_{i=1}^M$. We remark that one can define an observable $\pi_i$ to extract the $i$th component of the state $\boldsymbol{\phi}$ at time $t$ (that is, $\pi_i\l(\boldsymbol{\phi}\l(t; \boldsymbol{\phi}_0\r)\r)) = \boldsymbol{\phi}_i\l(t; \boldsymbol{\phi}_0\r)$ \cite{linDatadrivenModelReduction2021}). In general, observables $g$ must be square integrable functions of the full state $\boldsymbol{\phi}$ against some distribution; see below section \ref{sec:distribution}. We will use $\mathbf{g}$ (resp.~$\mathbf{g}\l(\boldsymbol{\phi}\r)$) to denote a vectorized observable $[g_1, g_2, \ldots g_M]^T$ (resp.~observations, $[g_1\l(\boldsymbol{\phi}\r), g_2\l(\boldsymbol{\phi}\r), \ldots g_M\l(\boldsymbol{\phi}\r)]^T$). 

\subsection{Discrete-time snapshot data} \label{sec:dataMatrix} Throughout this study, we assume that the full system has been simulated and observed at discrete times, and the observations form a data set for learning the Mori--Zwanzig operators. Specifically, we rely on samples drawn from a chosen initial distribution $\mu$ (see below section \ref{sec:distribution}) to collect the statistics of the unresolved information. The full system is prepared at $N \gg 1$ independently and identically sampled $\boldsymbol{\phi}_0^{[i]} \sim \mu$, $i=1\ldots N$. The system is then simulated and the values of the observations are registered at uniformly distributed discrete times $k\Delta$, $k=0\ldots K-1$. We refer to these observations as the \emph{snapshots}. Without loss of generality, we assume $\Delta=1$ by choosing an appropriate unit of time, unless it is specified otherwise. The full discrete-time observations thus form an $N \times M \times K$ data matrix $\bf D$ whose $(i,j,k)$-entry is $g_j \l(\boldsymbol{\phi}\l(k-1; \boldsymbol{\phi}_0^{[i]}\r)\r)$. For those systems which have natural or physical invariant measure, it is desirable to bypass sampling from a pre-specified initial distribution and use the long-time statistics instead. In this case, we generate a single but long trajectory from a randomly selected initial condition $\boldsymbol{\phi}_0^\text{r}$. After the initial transient time $t_c$, we make $K+N$ observations for every $\Delta=1$. We use these $K+N$ observations to inform the statistics of the invariant measure. The assumptions underlying the sufficiency of using $K+N$ observations are described below. The data matrix $\bf D$ in this case is defined as an $N\times M \times K$ matrix whose $(i,j,k)$-entry is $g_j \l(\boldsymbol{\phi}\l(t_c +i+k-1 ; \boldsymbol{\phi}_0^\text{r}\r)\r)$.

\subsection{\bf Discrete-time dynamical systems} To match to the discrete-time snapshot data, we   re-write the evolutionary equation \eqref{eq:ODE} into a discrete-time form:
\begin{equation}
    \boldsymbol{\phi}\l(t+1\r) = \mathbf{F}\l(\boldsymbol{\phi}(t)\r), \label{eq:disDyn}
\end{equation}
where $\mathbf{F}:\mathbb{R}^D \rightarrow \mathbb{R}^D$ is the flow map. We assume that the map exists and is unique for any finite time $t\in \mathbb{Z}_{\ge 0}$. Our formulation is readily applicable to those generically discrete-time dynamical systems, in which we begin with the discrete-time formulation Eq.~\eqref{eq:disDyn} instead of Eq.~\eqref{eq:ODE}.

\subsection{A chosen distribution and square integrable functions} \label{sec:distribution}We assume that the number of samples $N$ is sufficiently large to capture the statistics of the time-dependent distribution induced by the deterministic dynamics from the initial distribution $\mu$, i.e., the solution of the generalized Liouville equation \cite{gerlich1973VerallgemeinerteLiouvilleGleichung} corresponding to Eq.~\eqref{eq:ODE} :
\begin{equation}
    \frac{\partial}{\partial t} \rho\l(\boldsymbol{\phi}, t\r) = - \sum_{i=1}^N \frac{\partial}{\partial \boldsymbol{\phi}_i}\l[\mathbf{R}_i\l(\boldsymbol{\phi}\r)  \rho\l(\boldsymbol{\phi}, t\r)\r], \text{ with the initial data } \rho\l(\boldsymbol{\phi}, 0\r)=\mu\l(\boldsymbol{\phi}\r), \label{eq:Liouville}
\end{equation}
where $\mu\l(\boldsymbol{\phi}\r)$ is the probability density function of the measure $\mu$, assuming the measure $\mu$ is dominated by the Lebesgue measure in the state space $\mathbb{R}^D$ so $\mu\l(\dd \boldsymbol{\phi} \r) = \mu \l( \boldsymbol{\phi} \r) \dd \boldsymbol{\phi}$. 
When using the single trajectory to sample the ergodic distribution, we assume that $N$ is sufficiently large such that the collected $N+K$ samples faithfully capture the statistics of the ergodic measure, i.e., the stationary solution $\rho_\ast$ of Eq.~\eqref{eq:Liouville} satisfying \begin{equation}
    0 = \sum_{i=1}^N \frac{\partial}{\partial \boldsymbol{\phi}_i} \l[ \mathbf{R}_i\l(\boldsymbol{\phi}\r)  \rho_\ast\l(\boldsymbol{\phi}\r) \r]. \label{eq:LiouvilleStationary}
\end{equation} 
In both cases, we assume that each of the observables, $g_i$, is a square integrable function against the induced measure. That is, if one chooses to learn from a time-dependent distribution (Eq.~\eqref{eq:Liouville}), $\int_{\mathbb{R}^D} \rho\l(\boldsymbol{\phi}, t\r) g_i^2 \l(\boldsymbol{\phi}\r)  \dd \boldsymbol{\phi} < \infty$ for $0 \le t \le \l(K-1\r)\Delta $, or if one chooses to learn from a stationary distribution (Eq.~\eqref{eq:LiouvilleStationary}), $\int_{\mathbb{R}^D} \rho_\ast\l(\boldsymbol{\phi}\r)  g_i^2 \l(\boldsymbol{\phi}\r) $\\
$ \dd \boldsymbol{\phi} < \infty$.

\subsection{The projection operator} Suppose all the $M$ resolved observables are independent, we shall need another $D-M$ independent unresolved observables to uniquely specify the system's state. Denote these unresolved observables, which are also real-valued functions of the system's state $\boldsymbol{\phi}$, by $u_j:\mathbf{R}^D \rightarrow \mathbb{R}$, $j=1\ldots D-M$. Importantly, observations of these unresolved variables, $u_j \l(\boldsymbol{\phi}\r)$, are not accessible. Consider a real-valued function of the full state expressed in both the resolved and unresolved observations, $h\l(\l\{g_i\l(\boldsymbol{\phi}\r)\r\}_{i=1}^M, \l\{u_j\l(\boldsymbol{\phi}\r)\r\}_{j=1}^{D-M}\r)$.  A projection operator $\mathcal{P}$ is an operator mapping $h$ to another real-valued function $\l(\mathcal{P}h\r):\mathbb{R}^D \rightarrow \mathbb{R}$ which only depends on the resolved observations. Mori--Zwanzig formalism utilizes this projection operator and decomposes any function $h\l(\l\{g_i\l(\boldsymbol{\phi}\r)\r\}_{i=1}^M, \l\{u_j\l(\boldsymbol{\phi}\r)\r\}_{j=1}^{D-M}\r)$ into $\l(\mathcal{P}h\r)\l(\l\{g_i\l(\boldsymbol{\phi}\r)\r\}_{i=1}^M\r)+h_\perp\l(\l\{g_i\l(\boldsymbol{\phi}\r)\r\}_{i=1}^M, \l\{u_j\l(\boldsymbol{\phi}\r)\r\}_{j=1}^{D-M}\r)$, where $h_\perp := h-\mathcal{P}h$. Note that the domain of $\l(\mathcal{P}h\r)$ is still both the resolved and unresolved observations, except that the function does not explicitly depend on the unresolved $u_j\l(\boldsymbol{\phi}\r)$. A projection operator $\mathcal{P}$ must satisfy the property $\mathcal{P}^2 = \mathcal{P}$, so that for any $h$, $\mathcal{P}^n h = \mathcal{P} h$, $n\in \mathbb{N}$. Such a property is critical in the derivation of the MZ formalism \cite{zwanzig2001nonequilibrium}.

\subsection{Discrete-time Mori--Zwanzig formalism} Following the derivations in \cite{darveComputingGeneralizedLangevin2009,linDatadrivenModelReduction2021,lin2021DataDrivenLearningMori,tian2021DatadrivenLearningMori}, the discrete-time Mori--Zwanzig formalism prescribes the exact evolutionary equations of the observations, that for any initial condition $\boldsymbol{\phi}_0 \in \mathbb{R}^D$:
\begin{equation}
    {\bf g}_{n+1} \l(\boldsymbol{\phi}_0\r) ={\boldsymbol{ \Omega}}^{(0)} \l({\bf g}_n \l(\boldsymbol{\phi}_0\r)\r)+ \sum_{\ell=1}^n {\boldsymbol{ \Omega}}^{(\ell)} \l({\bf g}_{n-\ell} \l(\boldsymbol{\phi}_0\r)\r)  + {\bf W}_n \l(\boldsymbol{\phi}_0\r).  \label{eq:GLE}
\end{equation}
Here, ${\bf g}_{n}:\mathbb{R}^D \rightarrow \mathbb{R}^M $ is an $M\times 1$ vector function of the initial state $\boldsymbol{\phi}_0$ such that
\begin{equation}
    {\bf g}_n \l(\boldsymbol{\phi}_0 \r) \triangleq {\bf g} \l(\boldsymbol{\phi} \l( n; \boldsymbol{\phi}_0 \r) \r) \equiv \mathbf{g} \l(\mathbf{F}^n \l(\boldsymbol{\phi}_0\r)\r), \label{eq:forwardBackwardEqu}
\end{equation}
noting that $\mathbf{g}_0 \equiv \mathbf{g}$.
In other words, denote the finite-time ($\Delta$) Koopman operator \cite{koopma32dynamical} by $\mathcal{K}$, ${\bf g}_n := \mathcal{K}^n {\bf g}$. For clarity, we write each component of this vector function as  $\mathbf{g}_{n,i}:=\mathcal{K}^n g_i$ if needed. Note that $\mathbf{g}_{0,i}=\mathbf{g}_{i}$. Equation \eqref{eq:GLE}, referred to as the \emph{Generalized Langevin Equation} (GLE), is the central result of the Mori--Zwanzig formalism. It states that the vectorized observation at time $n+1$, ${\bf g} \l(\boldsymbol{\phi} \l( n+1; \boldsymbol{\phi}_0 \r) \r)$, can be decomposed into three parts: (1) a Markovian function $\boldsymbol{\Omega}^{(0)}:\mathbb{R}^M \rightarrow \mathbb{R}^M$ mapping the observations at time $n$ to an $M \times 1$ vector, (2) a series of functions $\boldsymbol{\Omega}^{(\ell)}:\mathbb{R}^M \rightarrow \mathbb{R}^M$, $\ell=1,2,\ldots$, mapping past observations with a lag $\ell$ to $M \times 1$ vectors, and (3) the \emph{orthogonal dynamics} that consist of observables ${\bf W}_n: \mathbb{R}^D \rightarrow \mathbb{R}^M $, $n=0, 1 \ldots$, mapping any initial state $\boldsymbol{\phi}_0 $ to an $M\times 1$ vector. Interested reader can find a detailed discussion and physical interpretation of the GLE \eqref{eq:GLE} in \cite{lin2021DataDrivenLearningMori}.

Because the GLE \eqref{eq:GLE} holds true for every initial condition $\boldsymbol{\phi}_0$, we can drop the function evaluation at $\boldsymbol{\phi}_0$ part and write succinctly an equation for the observables
\begin{equation}
    {\bf g}_{n+1} =\sum_{\ell=0}^n {\boldsymbol{ \Omega}}^{(\ell)} {\bf g}_{n-\ell}  + {\bf W}_n .  \label{eq:GLEo}
\end{equation}
In writing so, we overload $\boldsymbol{\Omega}^{(\ell)}$, which are functions in Eq.~\eqref{eq:GLE}, to denote operators in Eq.~\eqref{eq:GLEo}. That is, in Eq.~\eqref{eq:GLEo}, $\boldsymbol{\Omega}^{(\ell)}$ are functional operators mapping an $M \times 1$ vector function, each of whose component is an $L^2\l(\mathbb{R}^D, \mu\r)$-function, to another $M \times 1$ vector function, each of whose component is an $L^2\l(\mathbb{R}^D, \mu\r)$-function again. We thus refer to $\boldsymbol{\Omega}^{(\ell)}$ as the ``MZ operators''. We will adopt the slightly abused notation that the same symbol $\mathbf{\Omega}^{(\ell)}$ is used to denote both the functions in equations of observations (such as Eq.~\eqref{eq:GLE}) and operators in equations of observables (such as Eq.~\eqref{eq:GLEo}). 

The operators $\boldsymbol{\Omega}^{\l(\ell\r)}$ and $\mathbf{W}_n$ depend on the choice of the projection operator $\mathcal{P}$, the choice of the vectorized observable $\mathbf{g}$, and the finite-time ($\Delta$) Koopman operator $\mathcal{K}$ \cite{darveComputingGeneralizedLangevin2009,linDatadrivenModelReduction2021}:
\begin{align}
    \boldsymbol{\Omega}^{(n)} :={}& \mathcal{P} \mathcal{K} \l[\l(1-\mathcal{P}\r)\mathcal{K} \r]^n, \label{eq:Omegas}\\
    {\bf W}_n :={}& \l[\l(1-\mathcal{P}\r) \mathcal{K} \r]^{n+1} \mathbf{g}, \label{eq:Ws} 
\end{align}
with $n\in \mathbb{Z}_{\ge 0}$. Because $\mathcal{P} \mathbf{W}_n=0$ (from the fact $\mathcal{P}^2 = \mathcal{P}$), $\mathbf{W}_n$ is often referred to as the orthogonal dynamics. From these definitions, it is straightforward to identify the following equation which relates the memory operators to the orthogonal dynamics
\begin{equation}
    \boldsymbol{\Omega}^{(n)} \mathbf{g} \equiv \mathcal{P} \mathcal{K} \mathbf{W}_{n-1}, \forall  n \in \mathbb{N}.
\end{equation}
The meaning of the above equation is more transparent when we apply the operators to any state $\boldsymbol{\phi}_0\in \mathbb{R}^D$:
\begin{equation}
    \boldsymbol{\Omega}^{(n)} \l(\mathbf{g}\l(\boldsymbol{\phi}_0\r)\r) \equiv \mathcal{P} \mathcal{K} \mathbf{W}_{n-1} \l(\boldsymbol{\phi}_0\r)\equiv \mathcal{P} \l(\l(\mathbf{W}_{n-1}\circ \mathbf{F}\r) \l(\boldsymbol{\phi}_0\r) \r), \forall n \in \mathbb{N}, \label{eq:GFD}
\end{equation}
in which we applied the definition of the Koopman operator $\mathcal{K}$ and used $\circ$ to denote function composition. The above equation is referred to as the \emph{generalized fluctuation-dissipation relation} (GFD). Below, we will demonstrate the GFD enables data-driven learning of the operators $\boldsymbol{\Omega}^{\l(\ell\r)}$'s. 

\subsection{\bf Regression}\label{sec:regression} Here, we provide a short description of the regression analysis. In regression analysis, the goal is to identify an optimal model $f$ as a function of the independent variables $\mathbf{x}$ for predicting the dependent variables $\mathbf{y}$. We first focus on a scalar variable $y$, which can be one of the many components of a vector $\mathbf{y}$. The model $f$ is defined as a family of functions $f(\mathbf{x};\boldsymbol{\theta})$ of the independent variables $\mathbf{x}$, parametrized by an array of parameters $\boldsymbol{\theta}$. The assumption of the (homoscedastic) regression model is that, in the most general form, there exists a random variable $\varepsilon$ such that the dependent variable $Y$ can be expressed as
\begin{equation}
    Y\l(\mathbf{x}\r) = f\l(\mathbf{x}; \boldsymbol{\theta}_\ast\r) + \varepsilon.
\end{equation}
Here, we use the standard capital notation to denote a random variable $Y$. 
Given a list of paired realizations $\l\{\mathbf{x}^{[i]}, y^{[i]}\r\}_{i=1}^N$, the best-fit parameter $\boldsymbol{\theta}_\ast$ minimizes the differences between the model prediction $f\l(\mathbf{x}; \boldsymbol{\theta}\r)$ and the dependent variable $Y$. The best-fit parameter $\boldsymbol{\theta}_\ast$ is solved by minimizing a negative log-likelihood function that depends on the noise distribution. For example, for normally distributed $\varepsilon$, the negative log-likelihood function is proportional to the mean squared error (MSE), so the best-fit parameter is the minimizer
\begin{equation}
    \boldsymbol{\theta}_\ast := \argmin_{\boldsymbol{\theta}} \frac{1}{N} \sum_{i=1}^N \l(y^{[i]}-f\l(\mathbf{x}^{[i]};\boldsymbol{\theta}\r)\r)^2.
\end{equation}
When the dependent variables are vectors, that is, $\mathbf{y}$ is $P>1$ dimensional, the cost function can be chosen as the sum of the MSE of $M$ independent regressions:
\begin{equation}
    \boldsymbol{\theta}_\ast := \argmin_{\boldsymbol{\theta}} \frac{1}{N} \sum_{i=1}^N \sum_{j=1}^P \l(y_j^{[i]}-f_j\l(\mathbf{x}^{[i]};\boldsymbol{\theta}\r)\r)^2. \label{eq:MSEcost}
\end{equation}
We remark that the above equation comes with three assumptions about the noise model. First, in each component $j$ of the $P$-dimensional $\mathbf{y}$, the noise $\varepsilon_j$ is normally distributed. Second, the noise $\varepsilon_j$ are independently and identically distributed in each of the dimension. Third, the noise is homoscedastic, which means that the variance does not depend on the independent variable $\mathbf{x}$. We remark that the negative log-likelihood function can take a different form if one chooses a different noise model (e.g., correlated or heteroscedastic noise). Nevertheless, our results below are general. Thus, in the analysis below, we express a general cost function by 
\begin{equation}
    C\l(\boldsymbol{\theta}; \l\{\mathbf{x}^{[i]}, \mathbf{y}^{[i]}\r\}_{i=1}^N\r), \label{eq:generalCost}
\end{equation}
so the best-fit parameters $\boldsymbol{\theta}_\ast$ are obtained by solving a general optimization problem:
\begin{equation}
    \boldsymbol{\theta}_\ast := \argmin_{\boldsymbol{\theta}}   C\l(\boldsymbol{\theta}; \l\{\mathbf{x}^{[i]}, \mathbf{y}^{[i]}\r\}_{i=1}^N\r). \label{eq:bestFitPars}
\end{equation}

\subsection{Regression-based projection operators} \label{sec:regressionAsProjection} Using the notations defined above, we treat the observations at the time $s\ge 0$ as the independent variables, and the evaluations of any real-valued function $h:\mathbb{R}^D\rightarrow \mathbb{R}$ of the full state $\boldsymbol{\phi}\l(t; \boldsymbol{\phi}_0\r)$ at a later time $t>s$ as the dependent variable. The samples of the dependent and the independent variables are generated from the fully-resolved simulation, whose initial condition $\boldsymbol{\phi}_0^{[i]}$ is sampled from the distribution $\mu$. Note that the time-dependent full state $\boldsymbol{\phi}\l(s; \boldsymbol{\phi}_0^{[i]}\r)$ and $\boldsymbol{\phi}\l(t; \boldsymbol{\phi}_0^{[i]}\r)$ are neither resolved nor accessed, but the observables $\mathbf{g}_s\left(\boldsymbol{\phi}_0^{[i]}\right):=\mathbf{g}\l(\boldsymbol{\phi}\l(s;\boldsymbol{\phi}_0^{[i]}\r)\r)$ and $h_t \left(\boldsymbol{\phi}_0^{[i]}\right):= h\l(\boldsymbol{\phi}\l(t;\boldsymbol{\phi}_0^{[i]}\r)\r)$ are measured and registered as we simulate the dynamics. We use the samples to regress $h_t \left(\boldsymbol{\phi}_0^{[i]}\right)$ on $\mathbf{g}_s\left(\boldsymbol{\phi}_0^{[i]}\right)$ using a family of functions $f\l(\mathbf{g}_s\left(\boldsymbol{\phi}_0^{[i]}\right);\boldsymbol{\theta}\r)$ by minimizing a chosen cost function:
\begin{equation}
    \boldsymbol{\theta}_\ast := \argmin_\theta  C\l(\boldsymbol{\theta}; \l\{\mathbf{g}_s\l(\boldsymbol{\phi}_0^{[i]}\r), h_t \left(\boldsymbol{\phi}_0^{[i]}\right)\r\}_{i=1}^N\r).
\end{equation}
The best-fit $f\l( \cdot ; \boldsymbol{\theta}_\ast\r): \mathbb{R}^M \rightarrow \mathbb{R}$ is now a function that depends on only the resolved initial observations at time $s$, $\mathbf{g}_s\l(\boldsymbol{\phi}_0\r)$. It is clear that another regression of $f( \mathbf{g}_s\l(\boldsymbol{\phi}_0\r); \boldsymbol{\theta}_\ast)$ on the independent variables $\mathbf{g}_s\l(\boldsymbol{\phi}_0\r)$ result in $f(\cdot; \boldsymbol{\theta}_\ast)$ again. As such, regression is a projection of $h_t$ to $f\l(\mathbf{g}_s\l(\boldsymbol{\phi}_0\r); \boldsymbol{\theta}_\ast\r)$:
\begin{equation}
    \l(\mathcal{P} h_t \r)\l(\mathbf{g}_s\l(\boldsymbol{\phi}_0\r) \r) = f\l(\mathbf{g}_s\l(\boldsymbol{\phi}_0\r); \boldsymbol{\theta}_\ast\r).  \label{eq:regression-basedProjection}
\end{equation}
In summary, with the choice of using regression for projection, applying $\mathcal{P}$ to any function $h_t$ of the full-system state $\boldsymbol{\phi}$ results in the best-fit parametric function $f\l(\cdot; \boldsymbol{\theta}_\ast\r)$ of the reduced-order observations $\mathbf{g}_s\l(\boldsymbol{\phi}_0\r)$. With a cost function such as Eq.~\eqref{eq:MSEcost}, a straightforward generalization can be made to those vector functions $\mathbf{h}_t:\mathbb{R}^D\rightarrow\mathbb{R}^P$, that $\l(\mathcal{P}\mathbf{h}_t\r)\l(\cdot\r) = \mathbf{f}\l(\cdot;\boldsymbol{\theta}_\ast\r)$, where $\mathbf{f}$ is a parametric family of $\mathbb{R}^M\rightarrow \mathbb{R}^P$ functions. We remark that the regression-based projection operator may no longer be the projection operators commonly adopted in existing Mori--Zwanzig models, e.g.~Mori's, finite-rank, and Zwanzig's projection operators.

\section{Theoretical Results} \label{sec:theoreticalResults}

\subsection{Learning Mori--Zwanzig operators with regression-based projection using snapshot data} \label{sec:extractingKernels} 
We now use the snapshot data matrix $\mathbf{D}$ to learn the MZ operators $\boldsymbol{\Omega}^{\l(n\r)}$ and orthogonal dynamics $\mathbf{W}_n$, $n \in \mathbb{Z}_{\ge 0}$. As will be seen below, this can be achieved by iteratively performing a regression problem as a projection, quantifying the residuals as the orthogonal dynamics, and utilizing the GFD to define the next-order regression problem. We begin with setting $n=0$, the GLE \eqref{eq:GLE} with the expression Eq.~\eqref{eq:Omegas} states
\begin{align}
    {\bf g}_1 \l(\boldsymbol{\phi}_0\r)  = {}& \mathcal{P} \mathcal{K}  \l({\bf g} \l(\boldsymbol{\phi}_0\r)\r) + {\bf W}_0 \l(\boldsymbol{\phi}_0\r) = \mathcal{P} \l[{\bf g}_1\l(\boldsymbol{\phi}_0\r)\r] + {\bf W}_0 \l(\boldsymbol{\phi}_0\r). \label{eq:firstOmega}
\end{align}
for any initial state $\boldsymbol{\phi}_0$. We used the definition of the Koopman operator $\mathcal{K} \l({\bf g} \l(\boldsymbol{\phi}_0\r)\r)= {\bf g}_1 \l(\boldsymbol{\phi}_0\r) := \mathbf{g} \l(\boldsymbol{\phi}\l(1; \boldsymbol{\phi}_0\r)\r)$, which are the observations made at one time step ahead of the initial condition $\boldsymbol{\phi}_0$. These observations can be accessed in the snapshot data matrix: they are simply $\mathbf{D}\l(\cdot, \cdot, 2\r)$. So can the  $\mathbf{g}\l(\boldsymbol{\phi}_0\r)$, which are $\mathbf{D}\l(\cdot, \cdot, 1\r)$. The regression-based $\mathcal{P}$ now projects ${\bf g}_1$ to an optimal function of $\mathbf{g}\l(\boldsymbol{\phi}_0\r)$, by regressing $\mathbf{D}\l(\cdot, \cdot, 2\r)$ on $\mathbf{D}\l(\cdot, \cdot, 1\r)$. For example, using the MSE cost function \eqref{eq:MSEcost}, we solve the following optimization problem for projection:
\begin{subequations}\label{eq:optimization1}
\begin{align}
    C^{(0)}\l(\boldsymbol{\theta};\mathbf{D}\r):={}& \frac{1}{N} \sum_{i=1}^N \sum_{j=1}^M \l[\mathbf{D}\l(i,j,2\r)-\mathbf{f}_j\l(\l\{\mathbf{D}\l(i,k,1\r)\r\}_{k=1}^M;\boldsymbol{\theta}\r)\r]^2,\label{eq:optimization1a}\\
    \boldsymbol{\theta}_\ast^{(0)} :={}& \argmin_\theta C^{(0)}\l(\boldsymbol{\theta};\mathbf{D}\r), \label{eq:optimization1b}\\
    \boldsymbol{\Omega}^{(0)}\l(\cdot\r) ={}& \mathbf{f} \l(\cdot ; \boldsymbol{\theta}_\ast^{(0)}\r),\label{eq:optimization1c}
\end{align}
\end{subequations}
where $\mathbf{f}\l(\cdot; \boldsymbol{\theta}\r):\mathbb{R}^M \rightarrow \mathbb{R}^M$ is a family of regressional functions defined by the user. Here, the superscript $(0)$ is used to denote the cost function ($C^{(0)}$) and the solution of the parameters ($\boldsymbol{\theta}^{(0)}$) of the $0^\text{th}$ order regression analysis. Once the best-fit model is solved, the orthogonal dynamics $\mathbf{W}_0$ is identified as the \emph{residual} of the regression problem from Eq.~\eqref{eq:firstOmega}. That is, for the $i\text{th}$ sample initial condition $\boldsymbol{\phi}_0^{[i]}$, the $j$th component of the orthogonal dynamics at time $0$ is:
\begin{align}
    \l(\mathbf{W}_0\r)_j \l(\boldsymbol{\phi}_0^{[i]}\r) ={}&  \l[{\bf g}_1\l(\boldsymbol{\phi}_0^{[i]}\r)-\mathcal{P}\l(\mathbf{g}_1\l(\boldsymbol{\phi}_0^{[i]}\r)\r)\r]_j\nonumber \\
    ={}&
    \mathbf{D}\l(i,j,2\r) -  \mathbf{f}_j \l(\l\{\mathbf{D}\l(i ,k,1\r)\r\}_{k=1}^M ; \boldsymbol{\theta}_\ast^{(0)}\r).
\end{align}

We now use mathematical induction to show that we can extract $\boldsymbol{\Omega}^{(n+1)}$ and samples of $\mathbf{W}_{n+1}$, given $\boldsymbol{\Omega}^{(\ell)}$, $\ell=0\ldots n$ and the snapshot data matrix $\mathbf{D}$. We first use GLE \eqref{eq:GLEo} to express  $\mathbf{W}_n$ in terms of $\boldsymbol{\Omega}^{(\ell)}$ and $\mathbf{g}_{n-\ell}$:
\begin{align}
    {\bf W}_n  ={}& {\bf g}_{n+1} - \sum_{\ell=0}^n {\boldsymbol{ \Omega}}^{(\ell)} \left({\bf g}_{n-\ell}\right) 
\end{align}
Then, we use the GFD, Eq.~\eqref{eq:GFD} to construct $\boldsymbol{\Omega}^{(n+1)}$:
\begin{align}
    \boldsymbol{\Omega}^{(n+1)} \l(\mathbf{g}\r) ={}& \mathcal{P}\mathcal{K} \mathbf{W}_n =  \mathcal{P} \l[\mathcal{K}\mathbf{g}_{n+1} - \mathcal{K} \sum_{\ell=0}^n {\boldsymbol{ \Omega}}^{(\ell)} \l({\bf g}_{n-\ell}\r) \r] \nonumber \\
    = {}& \mathcal{P} \l({\bf g}_{n+2} - \sum_{\ell=0}^n {\boldsymbol{ \Omega}}^{(\ell)} \l({\bf g}_{n-\ell+1} \r) \r).
\end{align}
We evaluate the above function on an initial state $\boldsymbol{\phi}_0$ and use the definition Eq.~\eqref{eq:forwardBackwardEqu} to obtain
\begin{align}
    \boldsymbol{\Omega}^{(n+1)} \l(\mathbf{g} \left(\boldsymbol{\phi}_0\right)\r) = \mathcal{P} \l({\bf g}\left( \boldsymbol{\phi}\left(n+2;\boldsymbol{\phi}_0\right)\right) - \sum_{\ell=0}^n {\boldsymbol{ \Omega}}^{(\ell)} \l({\bf g}\left( \boldsymbol{\phi}\left(n-\ell+1;\boldsymbol{\phi}_0\right)\right)\r) \r).
\end{align}
In the right-hand side of this expression, $\boldsymbol{\Omega}^{\l(\ell\r)}$, $\ell=0\ldots n$, is the given best-fit parametric functions $\mathbf{f}\l(\cdot ; \boldsymbol{\theta}_\ast^{(\ell)}\r)$. Note that for a particular initial sample $\boldsymbol{\phi}^{[i]}$, all the observations in the above equation are registered in the data matrix $\mathbf{D}$:
\begin{algorithm}[t]
    \caption{Extracting the regression-based Mori--Zwanzig operators. \textbf{Input:} A set of observables $\mathbf{g}:\mathbb{R}^D \rightarrow \mathbb{R}^M$; the data matrix $\mathbf{D}$ of the snapshots of the observations (see Sec.~\ref{sec:dataMatrix}); a family of regressional functions $\mathbf{f}\left(\cdot;\right):\mathbb{R}^M \rightarrow \mathbb{R}^M$ parametrized by fitting parameters $\boldsymbol{\theta}$; the cost function of the regression analysis, $C\left(\boldsymbol{\theta}; \left\{\mathbf{y}^{[i]},\mathbf{x}^{[i]}\right\}_{i=1}^N\right)$ (see Sec.~\ref{sec:regression}).
    \textbf{Output:} The Markov operator $\boldsymbol{\Omega}^{(0)}:\mathbb{R}^M \rightarrow \mathbb{R}^M$; memory operators $\boldsymbol{\Omega}^{(n)}:\mathbb{R}^M \rightarrow \mathbb{R}^M$, $n \in \left\{1, 2, \ldots K-2\right\}$; orthogonal dynamics $\mathbf{W}_n$, $n\in \left\{1, 2, \ldots K-3\right\}$.}
    \label{alg:1}
    \begin{algorithmic}
    \State $\mathbf{x}_j^{[i]} \leftarrow \mathbf{D}\left(i,j,1\right)$
    \State $\mathbf{y}_{j,0}^{[i]} \leftarrow \mathbf{D}\left(i,j,2\right)$
    \State $\boldsymbol{\theta}_\ast^{(0)} \leftarrow \argmin_{\boldsymbol{\theta}} C\left(\boldsymbol{\theta}; \left\{\mathbf{y}_{0}^{[i]},\mathbf{x}^{[i]}\right\}_{i=1}^N\right) $
    \State $\boldsymbol{\Omega}^{(0)}\left(\cdot \right) \leftarrow \mathbf{f}\left(\cdot; \boldsymbol{\theta}_\ast^{(0)}\right)$
    \For{$n$ in $\left\{1, \ldots K-2 \right\}$}
        \State $\mathbf{y}_{j,n}^{[i]} \leftarrow \mathbf{D}\left(i,j,n+2\right) - \sum_{\ell=0}^{n-1} \mathbf{f}_j \l(\l\{\mathbf{D}\l(i,k,n-\ell+1\r)\r\}_{k=1}^M;\boldsymbol{\theta}_\ast^{\l(\ell\r)} \r)$ 
        \State $\mathbf{W}_{n-1}^{[i]} \leftarrow \mathbf{y}_{n}^{[i]}$
        \State $\boldsymbol{\theta}_\ast^{(n)} \leftarrow \argmin_{\boldsymbol{\theta}} C\left(\boldsymbol{\theta}; \left\{\mathbf{y}_{n}^{[i]},\mathbf{x}^{[i]}\right\}_{i=1}^N\right) $
        \State $\boldsymbol{\Omega}^{(n)}\left(\cdot \right) \leftarrow \mathbf{f}\left(\cdot; \boldsymbol{\theta}_\ast^{(n)}\right)$    
    \EndFor
    \State Output $\boldsymbol{\Omega}^{(n)}$ and $\mathbf{W}_n$
    \end{algorithmic}
    \end{algorithm}
\begin{subequations}
\begin{align}
    \mathbf{g}\left(\boldsymbol{\phi}_0^{[i]}\right) ={}&  \mathbf{D}\left(i,\cdot,1\right),\\ \mathbf{g}\left(\boldsymbol{\phi}\left(n+2;\boldsymbol{\phi}_0^{[i]}\right)\right) ={}& \mathbf{D}\left(i,\cdot,n+3\right),\\
    \mathbf{g}\left(\boldsymbol{\phi}\left(n-\ell+1;\boldsymbol{\phi}_0^{[i]}\right)\right) ={}& \mathbf{D}\left(i,\cdot,n-\ell+2 \right).
\end{align}
\end{subequations} 
Similar to Eq.~\eqref{eq:optimization1}, we now apply the projection operator $\mathcal{P}$, which is a regression on the corresponding variables:
\begin{subequations}\label{eq:optimization2}
\begin{align}
    \mathbf{y}^{[i]}_{j,n+1} :={}&  \mathbf{D}\l(i,j,n+3\r)-\sum_{\ell=0}^n \mathbf{f}_j \l(\l\{\mathbf{D}\l(i,k,n-\ell+2\r)\r\}_{k=1}^M;\boldsymbol{\theta}_\ast^{\l(\ell\r)} \r), \label{eq:nthOrderY} \\
    C^{(n+1)}\l(\boldsymbol{\theta};\mathbf{D}\r):={}& \frac{1}{N} \sum_{i=1}^N \sum_{j=1}^M \l[ \mathbf{y}^{[i]}_{j,n+1} - \mathbf{f}_j\l(\l\{\mathbf{D}\l(i,k ,1\r)\r\}_{k=1}^M;\boldsymbol{\theta}\r)\r]^2,\label{eq:optimizationG1} \\
    \boldsymbol{\theta}_\ast^{(n+1)} :={}& \argmin_{\boldsymbol{\theta}} C^{(n+1)}\l(\boldsymbol{\theta};\mathbf{D}\r), \label{eq:optimizationG2}\\
    \boldsymbol{\Omega}^{(n+1)} \equiv {}& \mathbf{f} \l(\cdot ; \boldsymbol{\theta}_\ast^{(n+1)}\r). \label{eq:optimizationG3}
\end{align}
\end{subequations}
Again, the superscript $(n+1)$ here to denotes the cost function ($C^{(n+1)}$), best-fit parameters $\boldsymbol{\theta}_\ast^{(n+1)}$, and best-fit parametric functions $\mathbf{f} \l(\cdot ; \boldsymbol{\theta}_\ast^{(n+1)}\r)$ of this $(n+1)$-th order regression analysis. By induction, the above procedure iteratively extracts $\boldsymbol{\Omega}^{(n)}\left(\cdot\right)=\mathbf{f} \l(\cdot ; \boldsymbol{\theta}_\ast^{(n)}\r)$, $n = 0 \ldots K-2$ (in total, $K-1$ operators), using the snapshot data matrix $\mathbf{D}$. Note that the family of regressional functions $\mathbf{f}$ is kept the same for all $\ell$ to ensure a consistent projection operator $\mathcal{P}$. Notably, both the orthogonality $\mathcal{P} \boldsymbol{W}_n =0$ and the GFD $\boldsymbol{\Omega}^{(n+1)}\left({\bf g}\right) = \mathcal{P}\mathcal{K} \mathbf{W}_n$ are built in by construction. We summarize the procedure in Algorithm \ref{alg:1}.

The above procedure iteratively learns the operators $\boldsymbol{\Omega}^{\l(\ell\r)}$ in the GLE \eqref{eq:GLEo}. The GLE is a particular form of memory-dependent dynamics, which is specific to the choices of the observables $\mathbf{g}$ and projection operator $\mathcal{P}$. There are other forms of memory-dependent dynamics which do not rely on an \emph{a priori} selected projection operator, for example, the Hankel-DMD method \cite{arbabiErgodicTheoryDynamic2017}. Numerical comparison of these methods and a detailed discussion will be provided below in Secs.~\ref{sec:KSmodel} and \ref{sec:discussion}. 

\subsection{Making predictions and modeling orthogonal dynamics} \label{sec:makingPrediction} Our goal is to iteratively use GLE \eqref{eq:GLE} to make multi-step predictions for the resolved observables, once the operators $\boldsymbol{\Omega}^{(\ell)}$'s are learned. For prediction, it is convenient to shift the time index and write the ``present time'' as $t=0$, assuming the system initiated from $\boldsymbol{\phi}_{-T}$ at $t=-T$ in the past ($T\in \mathbb{Z}_{\ge0}$). We assume that the present and a finite number of past times, $\mathbf{g}_{-\ell}\l(\boldsymbol{\phi}_{-T}\r)$, $\ell=0,1,\ldots H-1 \le T$, are given. We will make predictions for the resolved variables at a future horizon $t=k\ge 1$, $\mathbf{g}_{k}\l(\boldsymbol{\phi}_{-T}\r)$.
To achieve this, we first use a truncated GLE to predict $\mathbf{g}_1$:
\begin{equation}
    {\bf g}_{1}\l(\boldsymbol{\phi}_{-T}\r) = \sum_{\ell=0}^{H-1} {\boldsymbol{ \Omega}}^{(\ell)} \l({\bf g}_{-\ell}\l(\boldsymbol{\phi}_{-T}\r)\r)  + {\bf W}_T\l(\boldsymbol{\phi}_{-T}\r). \label{eq:GLEtruncated}
\end{equation}
This is a truncated GLE because only a finite length of the history is provided; terms ${\boldsymbol{ \Omega}}^{(\ell)} \l({\bf g}_{-\ell}\l(\boldsymbol{\phi}_{-T}\r)\r)$, $\ell=H \ldots T$ are neglected. We show below with numerical evidence that this truncation is justifiable.
Once $\mathbf{g}_1\l(\boldsymbol{\phi}_{-T}\r)$ is obtained, we accumulate the predicted $\mathbf{g}_1\l(\boldsymbol{\phi}_{-T}\r)$ into the series of observations, $\mathbf{g}_{-\ell}\l(\boldsymbol{\phi}_{-T}\r)$, $\ell=-1,0,1\ldots H-1 $, and use this accumulated series to make prediction of $\mathbf{g}_2\l(\boldsymbol{\phi}_{-T}\r)$ using the truncated GLE \eqref{eq:GLEtruncated} again:
\begin{equation}
    {\bf g}_{2}\l(\boldsymbol{\phi}_{-T}\r) = \sum_{\ell=0}^{H-1} {\boldsymbol{ \Omega}}^{(\ell)} \l({\bf g}_{1-\ell}\l(\boldsymbol{\phi}_{-T}\r)\r)  + {\bf W}_{T+1}\l(\boldsymbol{\phi}_{-T}\r).
\end{equation}
The procedure continues until the desired horizon $m\ge1$ is reached. Note that this procedure only needs the first $H$ operators $\boldsymbol{\Omega}^{(\ell)}$, $\ell=0\ldots H-1$. Clearly, the orthogonal dynamics $\mathbf{W}_n$, $n=T\ldots T+m-1$, are needed in the procedure. Recall that the orthogonal dynamics encode the information in the unresolved space (Eq.~\eqref{eq:Ws}), so they cannot be obtained by data-driven approaches. We must then rely on modeling. However, modeling the orthogonal dynamics is challenging: it is equivalent to modeling all the $D-M$ unresolved degrees of freedom. In addition, the modeling is specific to the full dynamical system, the choice of the resolved observables, and the choice of the projection operator. 

A best-case scenario is that the resolved observables are reasonably selected, and an expressive enough regression model is chosen. With a sufficiently long history length $H$, this leads to negligible orthogonal dynamics, i.e., the residuals of the regression analysis are small enough to be negligible. In the numerical experiments below, we assume this best-case scenario and proceed with a trivial zero-noise model $\mathbf{W}_n=0 $, $n \ge T$, despite that finite residuals are observed in the experiments. Prediction with $\mathbf{W}_n=0$, $n\ge T$ can be conceived as the conditional mean with respect to the distribution of under-resolved dynamics ($\mathbf{W}_n$, $n\ge T$). We remark that Price et al.~\cite{stinis21} also attacked the same problem using an operator algebraic expansion (of time $t$) to approximate the memory kernel, which is different from our purely data-driven approach. In addition, the GFD \eqref{eq:GFD} is guaranteed in our method, but not guaranteed in the $t$-expansion approximates.   

We also tested the standard way of fitting the residual $\mathbf{W}_n$ as independent and multivariate Gaussian distribution (independent between different time indices, and multivariate in the $M$ observables). The results are qualitatively similar but more noisy, because the prediction are based on randomly generated samples of the distribution. As the results with the Gaussian noise model are similar but less ``clean'', we will not present them in this article. 

\subsection{Linear and nonlinear regression-based projections} \label{sec:linearVnonlinear} 

We now elaborate a critical idea that distinguishes the subtle differences between two regression-based projections that appear nearly identical. We will refer to the first technique as ``linear regression-based projection'' and the second as ``nonlinear regression-based projection'', as we currently lack more precise terminology\footnote{We are open to reviewers' suggestion.}. Although the regression analysis for these two methods are highly similar, drastic differences exist in the procedures utilized to apply the trained models for making predictions. Consequently, the learned models can even produce qualitatively distinct long-time behaviors, as will be seen in the next Section.  

\begin{figure}[!t]
    \centering
    \includegraphics[width=0.94\textwidth]{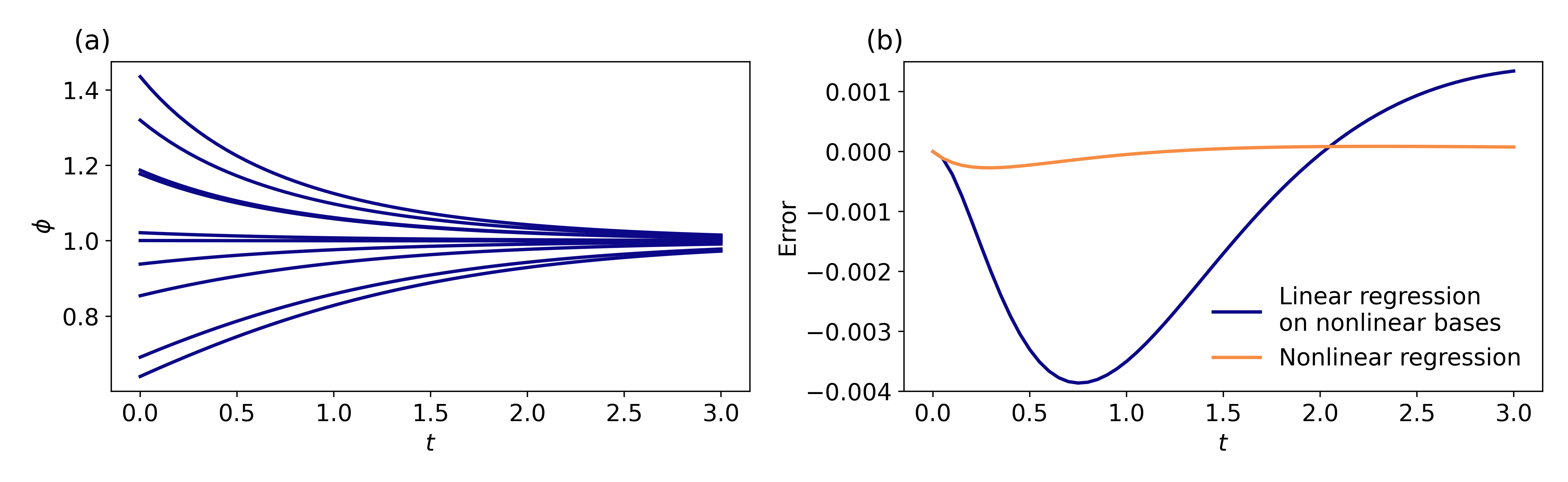}
    \caption{Learning the toy model $\dot{\phi}=-\phi^2 + \phi$. (a) Ten trajectories of the dynamics form randomly ($\sim 0.5+\text{Beta}(2,2)$) generated initial conditions. (b) The error of the prediction of the two regression-based projection operators.}
    \label{fig:phi1phi2}
\end{figure}

To illustrate this idea, we employ a one-dimensional toy model: $\dot{\phi}= - \phi^2 + \phi $ for some initial distribution $\phi(0) \sim \mu(\dd \phi)$. We chose $\mu =0.5+\text{Beta}(2,2)$ for this example. We generate $N=10^5$ trajectories in the time window $t \in [0, 3]$ and record snapshots every $\Delta=0.05$ (each trajectory contains $K=61$ snapshots, including the initial condition). Ten randomly selected trajectories are visualized in Fig.~\ref{fig:phi1phi2}(1). We will use the generated trajectories to perform two regression analyses, which we will later use for long-term predictions by the truncated GLE \eqref{eq:GLE}. For simplicity, we will only use the Markov term for making the predictions. This can be achieved by setting $\mathbf{W}_n=0$ for $n\in \mathbb{Z}{\ge0}=0$ during the prediction phase. Therefore, we will use the following equation recursively to make a multi-step prediction:
\begin{equation}
    {\bf g}_{m+1}(\boldsymbol{\phi}_0)={\boldsymbol{\Omega}}^{(0)}\left({\bf g}_m (\boldsymbol{\phi}_0)\right) \equiv \left [\left(\mathcal{P}\mathcal{K} \right)^m \mathbf{g} \right]\left(\phi_0\right). \label{eq:MarkovTransitionOnly}
\end{equation}
We would like to emphasize that the concept highlighted in this section is generalizable to systems that include Mori--Zwanzig memory. We focus on only the Markov transitions for the clarity of presentation. 

The first projection is in line with the approximate Koopman learning \cite{schmid2010DynamicModeDecomposition,williams2015DataDrivenApproximation} and data-driving learning of MZ operators with Mori's projection \cite{lin2021DataDrivenLearningMori}. We will refer to this method as the \textit{linear regression-based projection}. With this method, one first specifies a set of basis functions; we consider polynomial functions up to the quadratic order, ${\bf g}(\boldsymbol{\phi}) = \l [g_1(\phi), g_2(\phi), g_3(\phi) \r]^T$ with $g_1(\phi):= 1$, $g_2(\phi)=:\phi$, and $g_3(\phi)=:\phi^2$. The data matrix $\mathbf{D}$, whose entries are $D_{i,j,k}=\left[\phi \left(k-1 ; \phi_0^{[i]}\right)\right]^{j-1}$, is thus $10^5 \times 3 \times 61$. With this method, a regression on the three nonlinear basis functions is used as the projection operator. The linear regression model aims to identify the approximate Koopman operator $\boldsymbol{\kappa} \in \mathbb{R}^{3\times 3}$ that minimizes the mean squared error:
\begin{equation}
    \boldsymbol{\kappa}_\ast^{(0)} := \operatorname{argmin}_{\boldsymbol{\kappa}} \frac{1}{N} \sum_{i=1}^{N} \sum_{j=1}^3 \sum_{k=1}^{60}  \left (D_{i,j,k+1} -  \sum_{\ell=1}^3 \boldsymbol{\kappa}_{j,\ell} D_{i,\ell,k}\right)^2. \label{eq:optimizationKoopman}
\end{equation}
where $\boldsymbol{\kappa}_{j,\ell}$ is the $(j,\ell)$-component of $\boldsymbol{\kappa}$. The cost function involves all three functional bases, and there are $3\times 3$ parameters in the regression analysis. Note that the three functional bases are treated as \textit{independent variables}, although knowing $g_2(\phi):=\phi$ informs the others, for arbitrary $\phi$. Our numerical experiment showed that
\begin{equation}
    \boldsymbol{\kappa}_\ast^{(0)} = \left[\begin{array}{ccc} 
 +1.000 & +0.000 & -0.000  \\
 +0.002 & +1.044 & -0.046  \\
 -0.082 & +0.265 & +0.816  
 \end{array}\right]\label{eq:optimizedKoopman}
\end{equation}
With the learned $\boldsymbol{\kappa}_\ast$, we now make an $n$-step prediction into the future, provided a current state $\phi_0$. The nonlinear function evaluation of the initial state are $\mathbf{g}_0\left(\phi_0\right)= \left[\phi_0^0, \phi_0^1, \phi_0^2\right]^T$; here, the subscript of $\mathbf{g}$ represents the nonlinear functional bases are at time step $0$. Because $\mathcal{P} \mathcal{K} \mathbf{g} = \boldsymbol{\kappa}_\ast^{(0)} \cdot \mathbf{g}$ in this regression model, it is straightforward to propagate the functional bases $m$-steps, $m\ge 1$, into the future to make predictions:
\begin{equation}
    \mathbf{g}_m^\text{(pred)}\left(\phi_0\right) = \left(\left(\boldsymbol{\kappa}_\ast^{(0)}\right)^m \cdot \mathbf{g}_0 \right) \left(\phi_0\right). \label{eq:mStepPredictionKoopman}
\end{equation}
That is, directly applying $m$ times the linear transformation $\boldsymbol{\kappa}^\ast$ to the initial functional bases. Note that this method provides a closed dynamical system of the resolved observables ($g_0$, $g_1$, and $g_2$). Furthermore, the evolutionary equation of the resolved observables is linear. Thus, we can conceive the method as a \textit{linear closure scheme}. Such a procedure for making multi-step predictions is natural and commonly adopted in the Koopman framework \cite{williams2015DataDrivenApproximation,arbabi2018data,yeung2019learning,luschDeepLearningUniversal2018}. The procedure has a drawback: because the prediction of the evolution of the functional bases is linear, the nonlinear moment information is not preserved. For example, the predicted second moment of the system after $m>0$ steps is not necessarily the squared predicted first moment:
\begin{align}
  \l[{\bf g}_m^\text{(pred)}(\phi_0) \r]_3 ={}& \underbrace{\sum_{\ell=1}^3 \left(\left(\boldsymbol{\kappa}_\ast^{(0)} \right)^m \right)_{3,\ell} g_{\ell}\left(\phi_0\right)}_{\text{Predicted second moment}}  \\ \ne{}& \underbrace{\left[\sum_{\ell=1}^3 \left(\left(\boldsymbol{\kappa}_\ast^{(0)} \right)^m \right)_{2,\ell} g_{\ell}\left(\phi_0\right) \right]^2}_{\text{Squared predicted first moment}}=\l[{\bf g}_m^\text{(pred)}(\phi_0) \r]_2^2. \nonumber
\end{align}

The above equation Eq.~\eqref{eq:mStepPredictionKoopman} only uses the learned Markov kernel for making prediction. When we use $T$ Mori--Zwanzig memory terms to make prediction (see Sec.~\ref{sec:makingPrediction}):
\begin{subequations} \label{eq:mStepPredictionKoopmanm}
\begin{align}
    \mathbf{g}_{n+1}^\text{(pred)}\left(\boldsymbol{\phi}_{-T}\right) ={}& \sum_{\ell=0}^{T} \boldsymbol{\kappa}_\ast^{(\ell)}  \cdot \tilde{\mathbf{g}}_{n-\ell} \left(\boldsymbol{\phi}_{-T}\right),\, n=0,1,\ldots \\
    \tilde{\mathbf{g}}_{k} ={}& \left\{\begin{array}{ll}
    \mathbf{g}_k, & \text{if } k=-T, -T+1, \ldots -1, 0, \, \text{(Given history)} \\
    \mathbf{g}_k^\text{(pred)}, & \text{else.}\\
    \end{array}\right.
\end{align}
\end{subequations}
That is, the predicted values are recursively used as imputation when making the more than one step of predictions. Because of the linear nature of the projection operator, predictions based on Eqs.~\eqref{eq:mStepPredictionKoopmanm} are linear in $\left\{\mathbf{g}_k\right\}_{k=-T}^{0}$. Below, we shall refer to those predictions based on \eqref{eq:mStepPredictionKoopman} as \textit{linear projection without memory} and  \eqref{eq:mStepPredictionKoopmanm} as \textit{linear projection with memory}. 

The second method, which we term as the \textit{nonlinear regression-based projection}, is more commonly adopted in modeling. The idea is to use a regression model on a set of independent observables we care about. For those observables which can be determined by (possibly nonlinear) transformation of the independent observables, the explicit form of the transformations are used in the regression model. After learning, one predicts the evolution of the independent observables, using those explicit forms of transformation. With our toy example, there is only one independent observable ${\bf g} = \l[g(\phi)\r]^T=\l[\phi\r]^T$. As such, the data matrix $\mathbf{D}$, with entries $D_{i,1,k}=\phi \left(k-1; \phi_0^{[i]}\right)$, is $10^5 \times 1 \times 61$. To match the order of the first method, we consider a second-order polynomial regression as the projection operator. The family of parametric functions is $f_{\boldsymbol{\theta}}\left(\phi\right)=f_{\alpha,\beta,\gamma}\left(\phi\right) := \alpha + \beta \phi + \gamma \phi^2$ with three fitting parameters $\theta = \left(\alpha, \beta, \gamma\right)$. Note that the constant $1$ and the $\phi^2$ are expressed explicitly as transformation of $\phi.$ We aim to identify the optimal parameters minimizing the mean squared error of the one-step prediction
\begin{equation}
    \boldsymbol{\theta}^\ast := \operatorname{argmin}_{\boldsymbol{\theta}} \frac{1}{N} \sum_{i=1}^{N} \sum_{k=1}^{60} \left [D_{i,1,k+1} - f_{\boldsymbol{\theta}}\left(D_{i,1,k};\theta \right) \right]^2. \label{eq:optimizationNonlinear}
\end{equation}
The above optimization problem, Eq.~\eqref{eq:optimizationNonlinear}, is a sub-problem of Eq.~\eqref{eq:optimizationKoopman}. As such, the solution of Eq.~\eqref{eq:optimizationNonlinear} is identically the second row of the solution (Eq.~\eqref{eq:optimizedKoopman}): $\alpha^\ast=\kappa^{\ast}_{1,1}\approx 0.002$, $\beta^\ast=\kappa^{\ast}_{1,2}\approx 1.044$, and $\gamma^\ast=\kappa^{\ast}_{1,3}\approx -0.046$. Despite an identical optimization solution, the prediction based on this nonlinear regression model is distinct from the first method (Eq.~\eqref{eq:mStepPredictionKoopman}) because $\boldsymbol{\Omega}^{(0)}\left(g\right)  = \mathcal{P} \mathcal{K} g = \alpha_\ast^{(0)} + \beta_\ast^{(0)} g + \gamma_\ast^{(0)} g^2$ now:
\begin{equation}
    {\bf g}_m^\text{(pred)}\left(\phi_0\right) = \underbrace{f_{\boldsymbol{\theta}_\ast^{(0)}} \left(\ldots \left( f_{\boldsymbol{\theta}_\ast^{(0)}}\left({\bf g}_0\left(\boldsymbol{\phi}_0\right)\right)\right)\right)}_{m\text{ function compositions}}. \label{eq:mStepPredictionNonlinear}
\end{equation}
Such an recursive composition of a nonlinear mapping is more in line with most  modeling for nonlinear dynamical systems: we assume a structure of the vector field, and after fitting the free parameters, the best-fit model is treated as a nonlinear dynamical system (v.~a linear system of an augmented set of nonlinear functional bases, Eq.~\eqref{eq:mStepPredictionKoopman}). Nonlinear systems-identification methods, such as Sparse Identification of Nonlinear Dynamics (SINDy, \cite{Brunton16SINDy,Kaheman20SINDyPI,fasel2022ensemble}), also follow the recursive Eq.~\eqref{eq:mStepPredictionNonlinear} in prediction phase. Note that this method also provides a closed dynamical system of the resolved observables, but in this case there is only one resolved observable, $g\left(\phi\right)=\phi$, yielding a nonlinear evolutionary equation of $\phi$. Below, we refer to this method as a \textit{nonlinear closure scheme}. 

Analogous to Eq.~\eqref{eq:mStepPredictionKoopman}, Eq.~\eqref{eq:mStepPredictionNonlinear} only uses the learned Markov model for making prediction. To make prediction with $T$ Mori--Zwanzig memory terms:
\begin{subequations} \label{eq:mStepPredictionNonlinearm}
\begin{align}
    \mathbf{g}_{n+1}^\text{(pred)}\left(\boldsymbol{\phi}_{-T}\right) ={}& \sum_{\ell=0}^{T} \mathbf{f}_{\theta_\ast^{(\ell)}} \left( \tilde{\mathbf{g}}_{n-\ell} \left(\boldsymbol{\phi}_{-T}\right) \right),\, n=0,1,\ldots \\
    \tilde{\mathbf{g}}_{k} ={}& \left\{\begin{array}{ll}
    \mathbf{g}_k, & \text{if } k=-T, -T+1, \ldots -1, 0, \, \text{(Given history)} \\
    \mathbf{g}_k^\text{(pred)}, & \text{else.}\\
    \end{array}\right.
\end{align}
\end{subequations}
Similar to Eq.~\eqref{eq:mStepPredictionKoopmanm}, predicted values are recursively used as imputation when making the more than one step of predictions. With the linear nature of the projection operator, predictions based on Eqs.~\eqref{eq:mStepPredictionKoopmanm} are no longer linear in $\left\{\mathbf{g}_k\right\}_{k=-T}^{0}$. Below, we shall refer to those predictions based on \eqref{eq:mStepPredictionNonlinear} as \textit{nonlinear projection without memory} and  \eqref{eq:mStepPredictionNonlinearm} as \textit{nonlinear projection with memory}. 

To illustrate the difference between these two closure schemes, we used the learned $\boldsymbol{\kappa}_\ast^{(0)}$ and $\left(\alpha_\ast^{(0)}, \beta_\ast^{(0)}, \gamma_\ast^{(0)}\right)$ to predict the dynamics of the toy model every $\Delta=0.05$ up to $t=3.0$, using linear projection without memory (Eqs.~\eqref{eq:mStepPredictionKoopman}) and nonlinear projection without memory (\eqref{eq:mStepPredictionNonlinear}). In this demonstration, we set the initial condition of the test trajectory at $\phi(0)=1.4$. The error of the prediction to the true dynamics is illustrated in Fig.~\ref{fig:phi1phi2}(b), which shows that the nonlinear closure scheme performed better than the linear one. In the next section, we will see that predictions based on the linear closure scheme could even lose important qualitative features of the fully resolved dynamics.

The nonlinear projection is a better formulation for accommodating data-driven models with neural networks (NNs), even in the context of approximate Koopman learning \cite{schmid2010DynamicModeDecomposition,williams2015DataDrivenApproximation} and data-driving learning of MZ operators with Mori's projection \cite{lin2021DataDrivenLearningMori}. For these problems, the set of observables have to be specified \emph{a priori}. Because the identification of the optimal basis functions remains an open problem, several recently proposed methods leverage NNs for learning better observables \cite{liExtendedDynamicMode2017,Yeung2017LearningDN,luschDeepLearningUniversal2018,wehmeyerTimelaggedAutoencodersDeep2018}. However, NNs are nonlinear functions of their parameters, so the joint learning problem (i.e., simultaneously learning the functional basis and the approximate Koopman operator), often with additional nonlinear regularizers, is not a linear regression problem. Consequently, we cannot categorize the NN-based approach as a simple Mori's projection, for which our proposed algorithms \cite{lin2021DataDrivenLearningMori} can extract the MZ memory kernels. A projection defined by the nonlinear projection is the most natural and accurate mathematical framework to describe these NN-based learning methods. By identifying ``training a NN'' as a projection, the algorithm described in Sec.~\ref{sec:extractingKernels} prescribes a principled way to extract the memory kernels of these NN-based learning methods.

\subsection{\bf A special case: linear regression} We now show that the special choice of linear projection results in the discrete-time algorithm we proposed in \cite{lin2021DataDrivenLearningMori}. For any observation $\mathbf{g}\l(\boldsymbol{\phi}_0\r)$, linear regression postulates that the parametric form of $\mathbf{f}$ is a linear superposition of the independent variables:
\begin{equation}
    \mathbf{f}\l(\mathbf{g}\l(\boldsymbol{\phi}_0\r);\boldsymbol{\kappa}\r) := \boldsymbol{\kappa} \cdot \mathbf{g}\l(\boldsymbol{\phi}_0\r).
\end{equation}
Generally, we include the constant function that maps any $\boldsymbol{\phi}$ to $1$ in the vectorized observable $\mathbf{g}$ to accounts for the biased term in  linear regression analysis. Here, the parameters $\boldsymbol{\theta}$ is the list of entries of the $M\times M$ matrix $\boldsymbol{\kappa}$. The best fit of a linear regression is analytically tractable. Define the $N \times M$ data matrix $\mathbf{X}$ at time $k-1$ whose entries are $\mathbf{X}_{i,j}\l(k-1\r) := \mathbf{D}\l(i,j,k\r)$. The solution of Eqs.~\eqref{eq:optimization1} is 
\begin{align}
    \boldsymbol{\kappa}_\ast^{(0)} = \mathbf{C}\l(1\r) \cdot \mathbf{C}^{-1}\l(0\r), \label{eq:MoriOmega0}
\end{align}
where $\mathbf{C}\l(k\r) := \mathbf{X}^T{(k)}\cdot \mathbf{X}(0)$ is the empirical $k$-lag correlation matrix of the observations. Equation \eqref{eq:MoriOmega0} is precisely the extracted Markov term in \cite{lin2021DataDrivenLearningMori}, and the approximate Koopman operator \cite{schmid2010DynamicModeDecomposition,Schmid2011,williams2015DataDrivenApproximation}. As for the higher orders, Eq.~\eqref{eq:nthOrderY} can be expressed as
\begin{equation}
    \mathbf{y}_{n+1} := \mathbf{X}\l(n+2\r)-\sum_{\ell=0}^n \boldsymbol{\kappa}_\ast^{\l(\ell\r)} \cdot \mathbf{X}\l(n-\ell+1\r),
\end{equation}
leading to the solution of the minimization problem Eqs.~\eqref{eq:optimizationG1}-\eqref{eq:optimizationG3}:
\begin{equation}
    \boldsymbol{\beta}_\ast^{(n+1)} = \l[\mathbf{C}\l(n+2\r) - \sum_{\ell=0}^n \boldsymbol{\kappa}_\ast^{(\ell)} \cdot \mathbf{C}\l(n-\ell +1\r) \r] \cdot \mathbf{C}^{-1}\l(0\r), \label{eq:recursiveCorrelations}
\end{equation}
which is exactly the formula derived from the GLE with Mori's projection operator \cite{lin2021DataDrivenLearningMori}. Thus, the linear regression-based projection with memory is equivalent to Mori's and the finite-rank projection. In addition, the Markov transition kernel with the linear regression-based projection is the approximate Koopman operator learned by EDMD \cite{williams2015DataDrivenApproximation,lin2021DataDrivenLearningMori}. As such,  linear projection without memory is equivalent to the approximate Koopman analysis such as DMD \cite{schmid2010DynamicModeDecomposition} and EDMD \cite{williams2015DataDrivenApproximation}, and the linear projection with memory is equivalent to our proposed data-driven learning with Mori's projection \cite{lin2021DataDrivenLearningMori}. 

\begin{figure}[!t]
    \centering
    \includegraphics[width=0.95\textwidth]{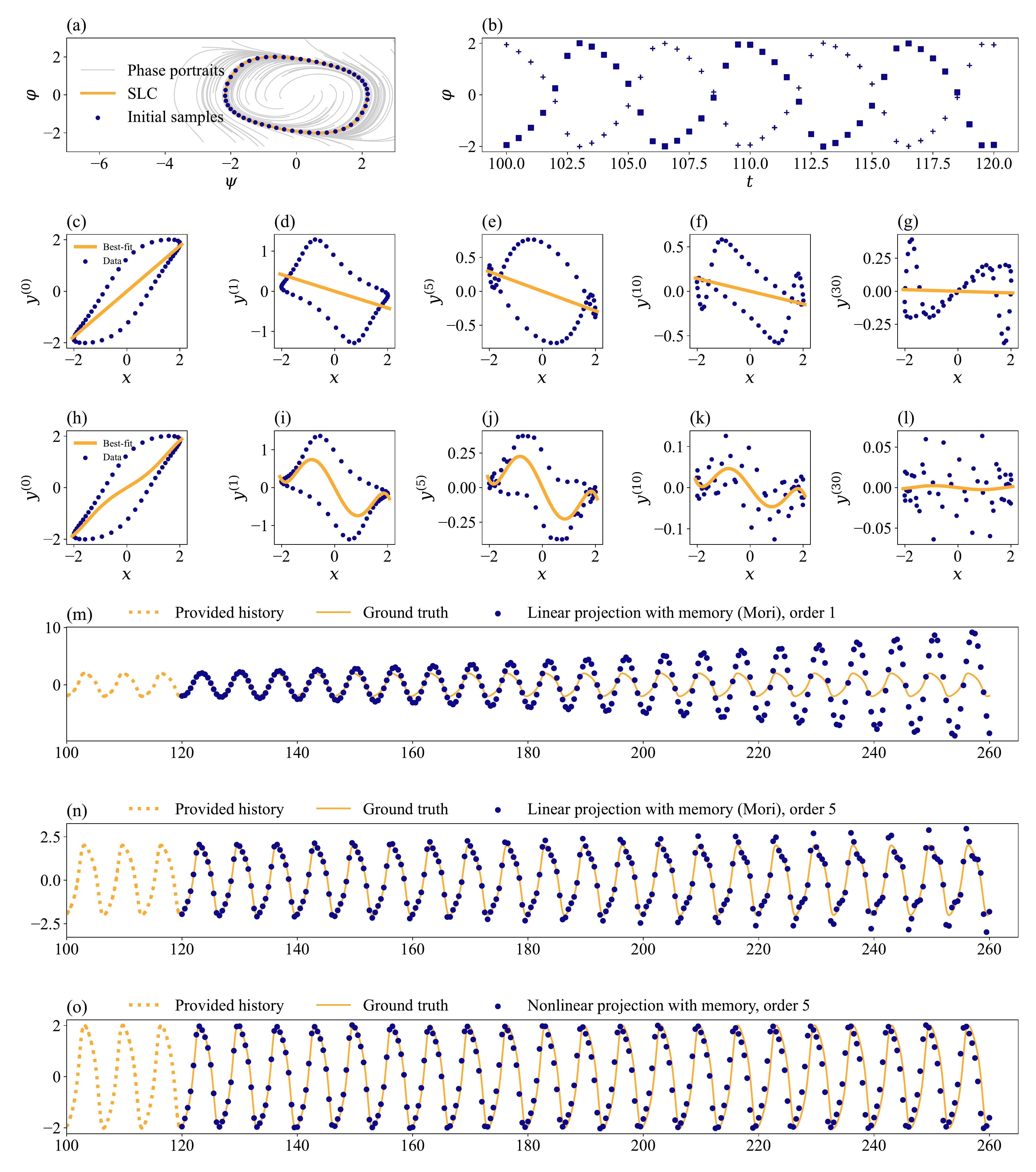}
    \caption{The Van der Pol oscillator. (a) The phase portraits of the full system, the stable limit cycle (SLC), and the initial samples of the data matrix $\mathbf{D}\l(\cdot, 1, 1\r)$; (b) Snapshot samples along two trajectories, $\mathbf{D}\l(1, 1, \cdot\r)$ and $\mathbf{D}\l(26, 1, \cdot\r)$; (c-g) Linear regression with a single observable, $f\l(x;\theta \r) := \theta x$; (h-l) Nonlinear polynomial regression with a single observable, $f\l(x;\boldsymbol{\theta}\r) := \sum_{i=0}^5 \theta_i x^i$; (m-o) Prediction by (m) linear projection with memory (Eq.~\eqref{eq:mStepPredictionKoopmanm}) on $\mathbf{g} \l(\boldsymbol{\phi}\r)\triangleq{} \l[1, \varphi\r]^T$;   (n) linear projection with memory (Eq.~\eqref{eq:mStepPredictionKoopmanm}) on $\mathbf{g} \l(\boldsymbol{\phi}\r)\triangleq{} \l[1, \varphi, \ldots \varphi^5\r]^T$;  and (o) nonlinear projection with memory (Eq.~\eqref{eq:mStepPredictionKoopmanm}) on $\mathbf{g} \l(\boldsymbol{\phi}\r)\triangleq{} \l[1, \varphi, \ldots \varphi^5\r]^T$;.}
    \label{fig:vanDerPol} 
\end{figure}

\section{Numerical Experiments} \label{sec:numerics}

\subsection{An illustrative example: the Van der Pol oscillator}

We now use the Van der Pol oscillator as an example to illustrate the regression-based learning of Mori--Zwanzig operators in action. The Van der Pol oscillator is a two-dimensional nonlinear system, which follows
\begin{subequations} \label{eq:VanDerPol}
\begin{align}
    \frac{\dd}{\dd t} \varphi(t) ={}& \mu \l(\varphi - \frac{\varphi^3}{3}\r) - \psi, \\
    \frac{\dd}{\dd t} \psi(t) ={}& \frac{1}{\mu} \varphi.
\end{align}
\end{subequations}
We chose $\mu=1$ in this illustrative example. The Van der Pol oscillator has a globally stable limit cycle, which is a one-dimensional manifold embedded in the $\boldsymbol{\phi}=\l(\varphi, \psi\r)$-state space, as shown in Fig.~\ref{fig:vanDerPol}(a). For illustrative purpose, we considered a single observable  $g\l(\boldsymbol{\phi}\r)\triangleq \varphi$ as the resolved dynamic variable. We used the \texttt{lsoda} integrator in \texttt{scipy.integrate.solve\_ivp} to numerically simulate the two-dimensional system from an arbitrarily chosen initial state $\boldsymbol{\phi}_0=(0,1)$ and a relaxation time $t_c=100$ for the system to relax to the limit cycle, whose period is roughly $6.663$. We chose $50$ initial conditions evenly sampled over one period; these initial conditions are also visualized in Fig.~\ref{fig:vanDerPol}(a). Then, we collected snapshot data set $\mathbf{D}$ by measuring $\varphi \l(t\r)$ from $t=0$ to $t=20$ every $\Delta=0.5$ from each of the initial conditions; two of the trajectories as shown in Fig.~\ref{fig:vanDerPol}(b). We used the mean-square error as the cost function for the regression models below. In Figure \ref{fig:vanDerPol}(c)-(g), we show how Mori's projection operator (i.e., linear projection) is applied to the snapshot data set $\mathbf{D}$. In Figure \ref{fig:vanDerPol}(c), the paired samples are simply $x=\mathbf{D}\l(\cdot,1,1\r)$ and $y=\mathbf{D}\l(\cdot,1,2\r)$. That is, the $x$'s are the observed $\varphi$ at the time $t=0$, and the $y^{(0)}$'s are the observed $\varphi$ one step ahead at the time $t=\Delta$. Because the system is not fully resolved, there can be two possible $y$'s for most of the measured $x$ on the stable limit cycle (see Fig.~\ref{fig:vanDerPol}(a)). Approximate Koopman and Mori's projector on a single observable is equivalent to a linear regression using $f(x;\theta)=\theta x$ to fit the data. The best-fit line is the projected model, $\theta_\ast^{(0)} x$. The slope of the best-fit model, $\theta_\ast^{(0)}$, is the approximate Koopman operator \cite{schmid2010DynamicModeDecomposition,Schmid2011,williams2015DataDrivenApproximation} and equivalently the $1\times 1$ Markov transition matrix in the context of data-driven Mori--Zwanzig formalism \cite{lin2021DataDrivenLearningMori}. After the first regression is made, the learned model $\theta_\ast^{(0)}$ is used to predict the system from $t=\Delta$ to $2\Delta$,  and the residual $ \mathbf{D}\l(\cdot,1,3\r) - \theta_\ast^{(0)} \mathbf{D}\l(\cdot,1,2\r)$ is assigned to $y^{(1)}$. The next linear regression on the independent variable, which remains as $x=\mathbf{D}\l(\cdot,1,1\r)$, is performed (Fig.~\ref{fig:vanDerPol}(d)), and the slope of the best-fit line is the first discrete-time $1\times 1$ memory kernel \cite{lin2021DataDrivenLearningMori}. The process repeats, and in Figs.~\ref{fig:vanDerPol}(e-g), we showed the $5$th, $10$th, and the $30$th linear regression. In Figs.~(h-l), we used fifth-order nonlinear polynomial regression instead. As expected, the more expressive polynomial regression led to smaller errors than that of the linear regression.

Figures (m-o) shows the predictions of the learned models. First, we simulated the full system \eqref{eq:VanDerPol} from a different initial condition $\boldsymbol{\phi}_0=(1,0)$ until $t=t_c=100$ for relaxing the system to the stable limit cycle. Then, we measure $\varphi(t+t_c;\boldsymbol{\phi}_0)$ every $\Delta=0.5$ until $t=260$ as the ground truth of the dynamics. We considered three learned models that are more expressive than a single linear observable ($\phi$), including and referred to as (1) linear projection with memory (order 1): two polynomial functions as the observables, $\mathbf{g} \l(\boldsymbol{\phi}\r)\triangleq{} \l[1, \varphi\r]^T$, predicted by Eq.~\eqref{eq:mStepPredictionKoopmanm}; (2) linear projection with memory (order 5): six lowest-order polynomial functions as the observables, $\mathbf{g} \l(\boldsymbol{\phi}\r)\triangleq{} \l[1, \varphi, \varphi^2, \ldots \varphi^5\r]^T$), predicted by Eq.~\eqref{eq:mStepPredictionKoopmanm}; (3) nonlinear projection with memory (order 5): six lowest-order polynomial functions as the observables, $\mathbf{g} \l(\boldsymbol{\phi}\r)\triangleq{} \l[1, \varphi, \varphi^2, \ldots \varphi^5\r]^T$), predicted by Eq.~\eqref{eq:mStepPredictionNonlinearm}. Given the snapshots in $t\in \l(100, 120\r)$, we used 
Eqs.~\eqref{eq:mStepPredictionKoopmanm} or Eqs.~\eqref{eq:mStepPredictionNonlinearm} to predict the dynamics at $t \ge 120$. Figures \ref{fig:vanDerPol}(m) and (n) shows that the error due to the negligence of the orthogonal dynamics accumulates and eventually destabilize the prediction in Mori's projection formalism. In contrast, Figure \ref{fig:vanDerPol}(o) shows that the nonlinear regression-based projection operator can reasonably approximate the true dynamics for a longer predictive horizon, using only the resolved snapshots.

We remark that this example was deliberately constructed as a minimal example. Our intention is to use a simple system and its transparent visualizations (Fig.~\ref{fig:vanDerPol}) to illustrate the iterative learning of the Mori--Zwanzig operators in action.

\subsubsection{Lorenz (1963) system}
Our second example is the Lorenz ((1963)) model \cite{lorenzDeterministicNonperiodicFlow1963}, which is a three-dimensional system
\begin{subequations}
\begin{align}
    \frac{\dd}{\dd t}\varphi\l(t\r) ={}& \sigma \l(\psi - \varphi\r), \\
    \frac{\dd}{\dd t}\psi\l(t\r) ={}& \varphi \l(\rho - \chi\r) - \psi  \\
    \frac{\dd}{\dd t}\chi\l(t\r) ={}& \varphi \psi - \beta \chi. 
\end{align}
\end{subequations}
We adopted the original model parameter values which E.~Lorenz chose, $\l(\sigma, \rho, \beta\r)=\l(10, 28, 8/3\r)$. Note that the full phase-space state $\boldsymbol{\phi}$ is $\left[\varphi, \psi, \chi\right]$. We again considered the ``reduced-order model'' as the $\varphi(t)$ alone. A long trajectory with an arbitrarily selected initial condition $\boldsymbol{\phi}_0^\text{tr}=\l(0.01, 1, 10\r)$ was generated using \texttt{scipy.integrate.} \texttt{solve\_ivp} for data-driven learning. We discarded the initial transient before $t=t_c=1000$, and collect $10^6$ snapshots of $\varphi$ every $\Delta=0.01$ until $t=t_c + 10^4$, as the samples on the strange attractor. A test trajectory was generated from a different initial condition $\boldsymbol{\phi}_0^\text{te}=\l(0, 1, 2\r)$, and also with $t_c=1000$. Again,  we used the mean-square error as the cost function for the regression models.

We considered five projection operators and compared their performances, ordered by the complexity (expressiveness) of the statistical model below. First, we again considered the order-1 linear projection operator, i.e., linear regression on the vectorized observable defined by $\mathbf{g}\l(\boldsymbol{\phi}\r)=\l[1, \varphi\r]^T$. We will refer to this as the Mori (order-1). Without the memory effects, such a projection operator corresponds to the DMD analysis \cite{schmid2010DynamicModeDecomposition,Schmid2011}. Next, we considered the order-5 linear projection operator but on five polynomial basis function, that is, linear regression on the vectorized observable which maps the full state $\boldsymbol{\phi}$ to $\l[1, \varphi, \varphi^2 \ldots \varphi^5\r]^T$. We will refer to this as the Mori (order-5). Without the memory, such a projection operator corresponds to an EDMD analysis \cite{williams2015DataDrivenApproximation}. The third projection operator we considered is an order-5 nonlinear projection defined as the fifth-order polynomial regression on $\varphi$. We will refer to this as the polynomial regression (order-5). Next, we considered a nonlinear projection operator defined by a cubic (degree $=3$) spline regression with ten knots (including boundary knots), evenly distributed between 1.5 times of the minimum and maximum of the collected $\varphi$ samples. After the spline features are generated, we performed a ridge regression with the $L^2$-regularization parameter set at $\lambda=10$. Finally, we considered a simple neural network as the nonlinear projection operator. The architecture of the neural networks is two fully-connected feed-forward layers, each of which contains five artificial neurons, and a third linear layer with only one neuron. In this study, we always chose the hyperbolic tangent function as the activation function of the artificial neurons. 

\newcommand{\FigResMemory}{{\bf \red{Figure 2} }}
\newcommand{\FigTrajectory}{{\bf \red{Figure 3} }}
\newcommand{\FigDistributions}{{\bf \red{Figure 3} }}
\newcommand{\FigErrors}{{\bf \red{Figure 4} }}
\newcommand{\FigErrorPaths}{{\bf \red{Figure 5} }}

\begin{figure}
    \centering
    \includegraphics[width=0.75\textwidth]{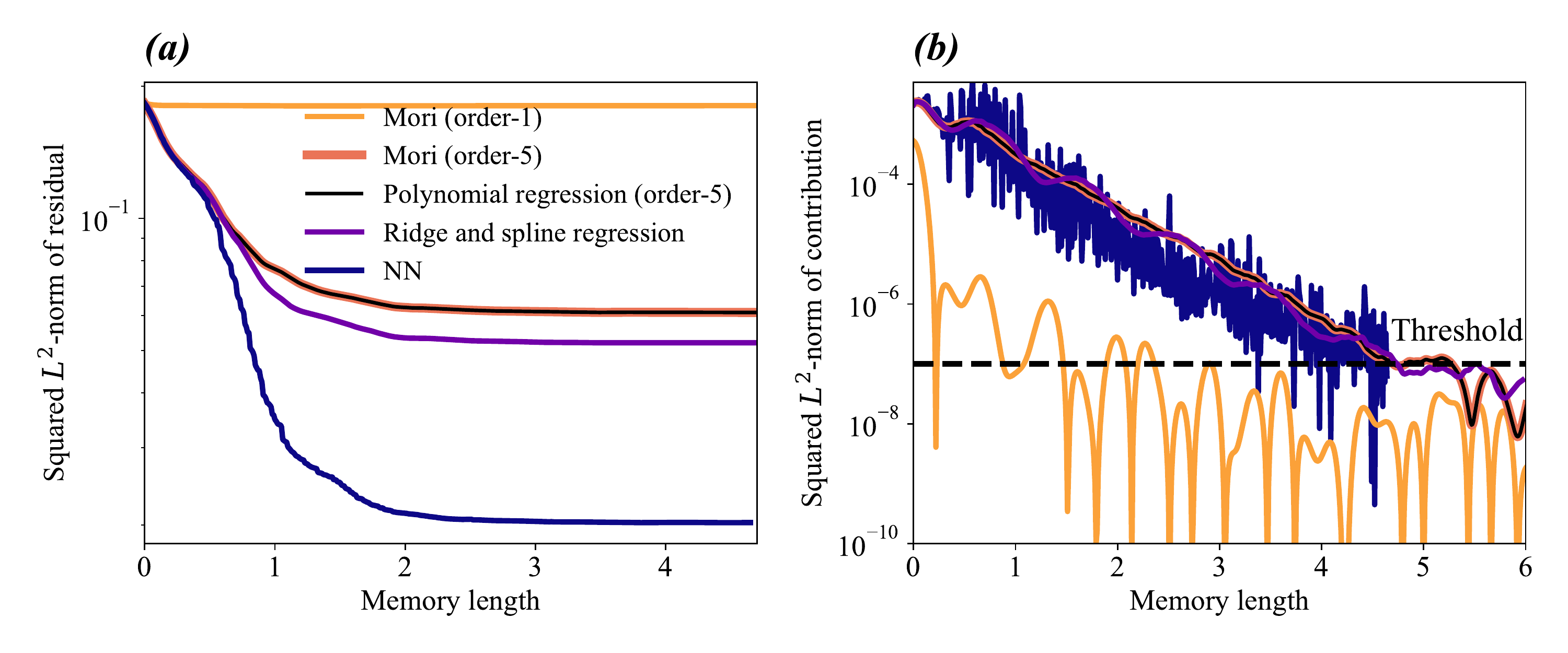}
    \caption{The Lorenz (1963) system. (a) The squared $L^2$-norm of the validation residual as the memory length increases. (b) The squared $L^2$-norm contribution at different memory length. The threshold is set at $10^{-7}$ for determining the memory length. Note that in these visualizations, Mori (5) and Nonlinear (5) are identical. The memory length is in the physical time unit, i.e., memory length $1$ corresponds to $100$ discrete-time steps ($\Delta=0.01$). }
    \label{fig:Lorenz63Res} 
\end{figure}

We used the measured snapshots to learn the operators $\boldsymbol{\Omega}^{\l(\ell\r)}$'s for each of the selected projection operators. Figure \ref{fig:Lorenz63Res}(a) shows the magnitude of the summary residual, defined as the averaged mean-squared error (MSE) of the best $\ell$th regression model on the test trajectory. Recall that $\ell$ is the index of the memory terms (cf.~Eq.~\eqref{eq:GLE}). Thus, we plotted the summary residual as a function of the physical time of the memory length ($\ell \Delta$) in Fig.~\ref{fig:Lorenz63Res}(a). In addition, for linear projections, the full optimization problem is the sum of the MSE of all the observables, in contrast to only the MSE of $\varphi$ in those nonlinear projections. To make a consistent comparison, we only quantified and reported the MSE of the component $\varphi$ for linear projections in the figure. It is observed that with the increased complexity of the regression model, the residual with memory contribution can decrease significantly, from order-1 linear projection ($\approx 2\times 10^{-1}$) to the expressive neural network ($\approx 2\times 10^{-2}$). Note that without the memory contribution ($\ell=0$), the magnitudes of the summary residual are indistinguishable; the significant improvement comes only with a finite-time memory contribution, $\approx 2.0$. This observation justified the advantage of the Mori--Zwanzig's memory-dependent formulation, especially for more expressive models.

\begin{figure}[!!t]
    \centering
    \includegraphics[width=\textwidth]{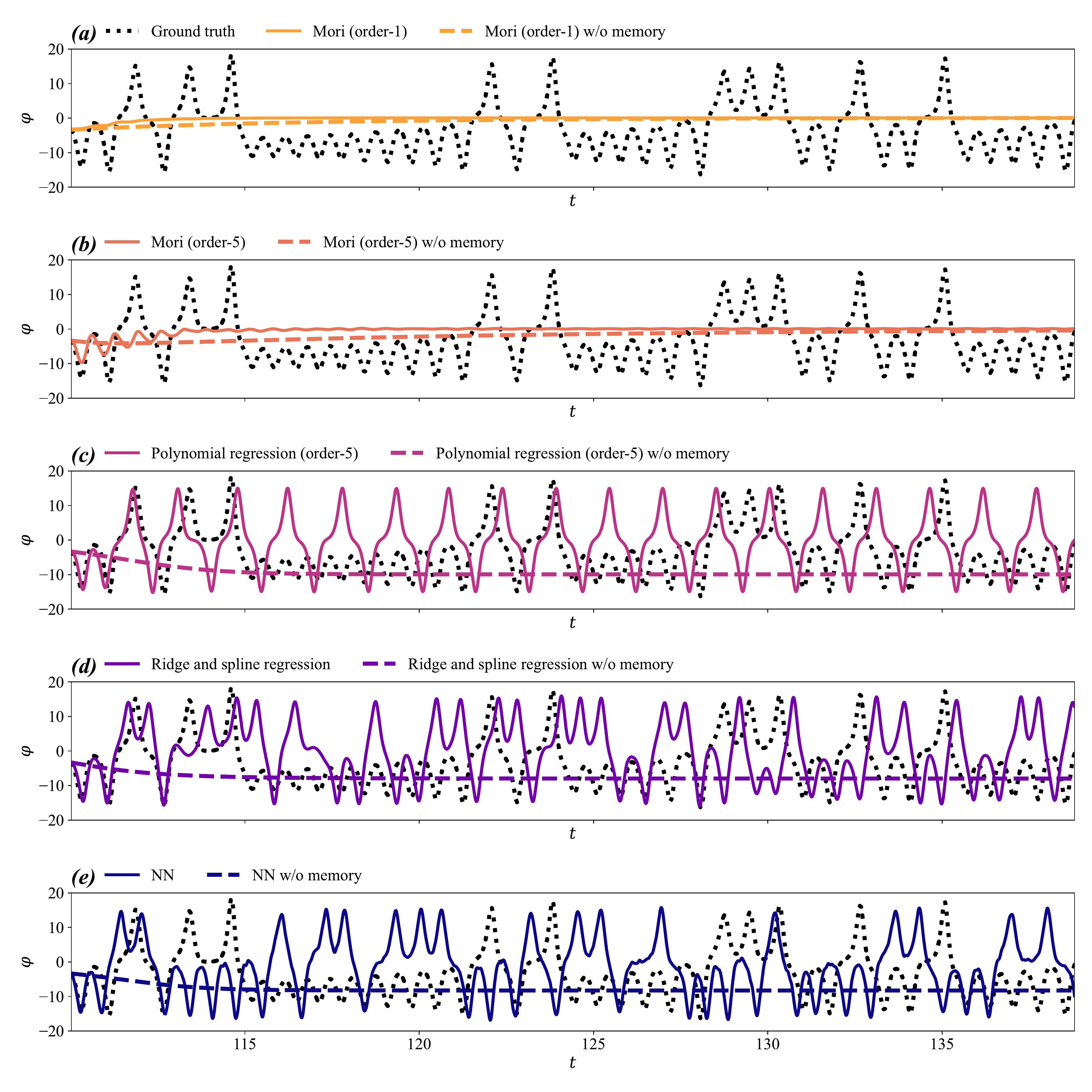}
    \caption{The Lorenz (1963) system. The predicted trajectories for the learned models based on (a) Mori (1) projection operator, (b) Mori (5) projection operator, (c) a projection by fifth-order polynomial regression on $\varphi$, (d) a projection by ridge and spline regression, and (e) a projection by a neural network (NN).}
    \label{fig:Lorenz63Trajs} 
\end{figure}

To investigate the contribution of the memory effects, we plotted the $L^2$-norm of the memory contribution in Eq.~\eqref{eq:GLE} ($\Omega^{(\ell)}\l(\boldsymbol{\phi}\l(k\Delta+t_c; \boldsymbol{\phi}_0^\text{te}\r)\r)$, $k=0,1,\ldots K+N$) averaged over the collected snapshots of the test trajectory. That is, we chose each snapshot along the long trajectory as a sample of the initial condition, and computed the memory contribution $\ell$ steps into the future. For order-1 linear projection operators, we only quantified the memory contributions to $\varphi$, i.e., $\l[\boldsymbol{\Omega}^{(\ell)} \cdot \boldsymbol{\phi} \l(k \Delta+t_c; \boldsymbol{\phi}_0^\text{te}\r)\r]_2$ (note that the first component is the constant function $1$). Figure \ref{fig:Lorenz63Res}(b) shows the memory contribution with various projection operators. We observed that memory contributions decayed as the memory index $\ell$ increased, regardless of the regression method. This observation suggests that a finite length ($H$) of the memory is sufficient, as the contribution of past history after a certain timescale would be small and negligible. As a result, in the following analysis, we chose a finite memory contribution by thresholding the squared $L^2$-norm of the memory contribution at $10^{-7}$. The corresponding memory lengths are: $H=235$ snapshots for order-1 linear projection; $H=469$ snapshots for order-5 linear projection; $H=469$ snapshots for the order-5 polynomial regression; $H=469$ snapshots for the spline regression; $H=469$ snapshots for neural network. We note that the number of discrete-time memories is much higher than what would have been needed for time-embedding techniques, such as Takens delay embedding \cite{takens81detecting}. This provides the numerical evidence that the Mori--Zwanzig memory is not the same memory contribution in Takens delay embedding technique, a point that we will elaborate in Sec.~\ref{sec:discussion}.

After the memory kernels ($\boldsymbol{\Omega}^{(\ell)}$) are learned, we turned our attention to the prediction of the learned models. We took the past history of the reduced-order variable $\varphi(1000 \le t\le 1050)$ and use Eq.~\eqref{eq:GLEtruncated} to make predictions for $t > 1050$, assuming the orthogonal dynamics $\mathbf{W}_n=0$. For comparison, we also made prediction from memory-less models, that is, setting $\boldsymbol{\Omega}^{(\ell \ge 1)}=0$. Note that the order-1 and order-5 linear without memory corresponds to DMD \cite{schmid2010DynamicModeDecomposition} and EDMD \cite{williams2015DataDrivenApproximation} respectively.

\begin{figure}[!t]
    \centering
    \includegraphics[width=\textwidth]{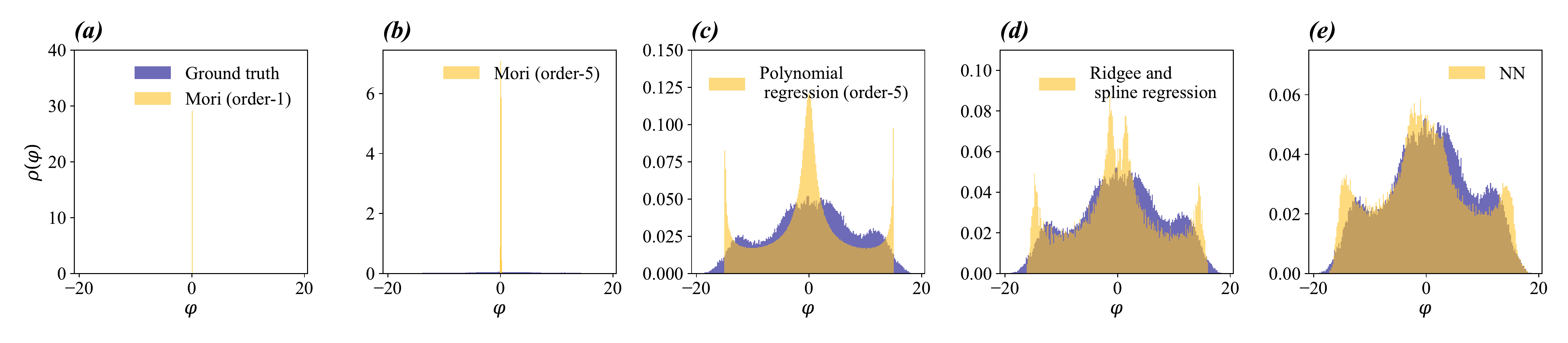}
    \caption{Long-time statistics of the MZ model predictions of the Lorenz (1963) model.}
    \label{fig:Lorenz63PDF} 
\end{figure}  

Figure \ref{fig:Lorenz63Trajs} shows the predicted trajectories by the learned models, using \ref{eq:mStepPredictionKoopmanm} for linear projections and \ref{eq:mStepPredictionNonlinearm} for nonlinear projections . For all the methods, we found that the memory-less models failed to capture the characteristics of the resolved dynamics other than the long-time mean. For order-1 and order-5 linear projections, a similar behavior of relaxation to the mean was observed. This should not be a surprise: as detailed in our recent work \cite{lin2021DataDrivenLearningMori}, linear projections without the unresolved orthogonal dynamics is functionally identical to using the two-time correlation functions as the propagators for prediction. The two-time correlation function $\mathbf{C}_{i,j}(t)$ decays to 0 because the chaotic dynamics de-correlates at a large lag $t \gg 1$. Consequently, linear projections also made predictions that converge to the mean behavior in the long-time limit. They were able to predict the mean because a constant function $1$, which can capture the constant bias, was included in the set of observables. Notably, linear projections with the memory kernels learned to predict a transient oscillation which crudely characterized the short-time dynamics of the Lorenz (1963) model. The drawback of the steady long-time prediction was significantly improved when we used nonlinear projection. The fifth-order polynomial regression showed a much better prediction in a short horizon, and an oscillation between the two wings of the Lorenz butterfly. Notably, the oscillation is not chaotic as it is in the fully resolved system. The prediction based on the spline regression, which is an intermediate model between the nonlinear polynomial regression and the NN-based regression models, began to exhibit chaotic behaviors in the long-time limit. Finally, the prediction of the one-dimensional reduced-order model based on an NN was qualitative similar to the chaotic dynamics of the fully-resolved Lorenz (1963) model.

To quantify the performance of the predictions, we first collected the long-time statistics by evolving Eq.~\eqref{eq:GLEtruncated} for 150,000 steps. We plotted the distributions of the collected trajectories based on the learned models in Fig.~\ref{fig:Lorenz63PDF}; such empirical distributions would be the stationary distribution if the process is ergodic. Clearly, linear projections showed an almost $\delta$-like distribution at $\varphi=0$ because their predictions converged to the mean in a finite timescale. All the predictions based on nonlinear projections reasonably approximated the empirical distribution of $\varphi$. We observed that the neural network out-performed the spline regression, which out-performed the plain fifth-order polynomial regression. 

\begin{figure}
    \centering
    \includegraphics[width=0.75 \textwidth]{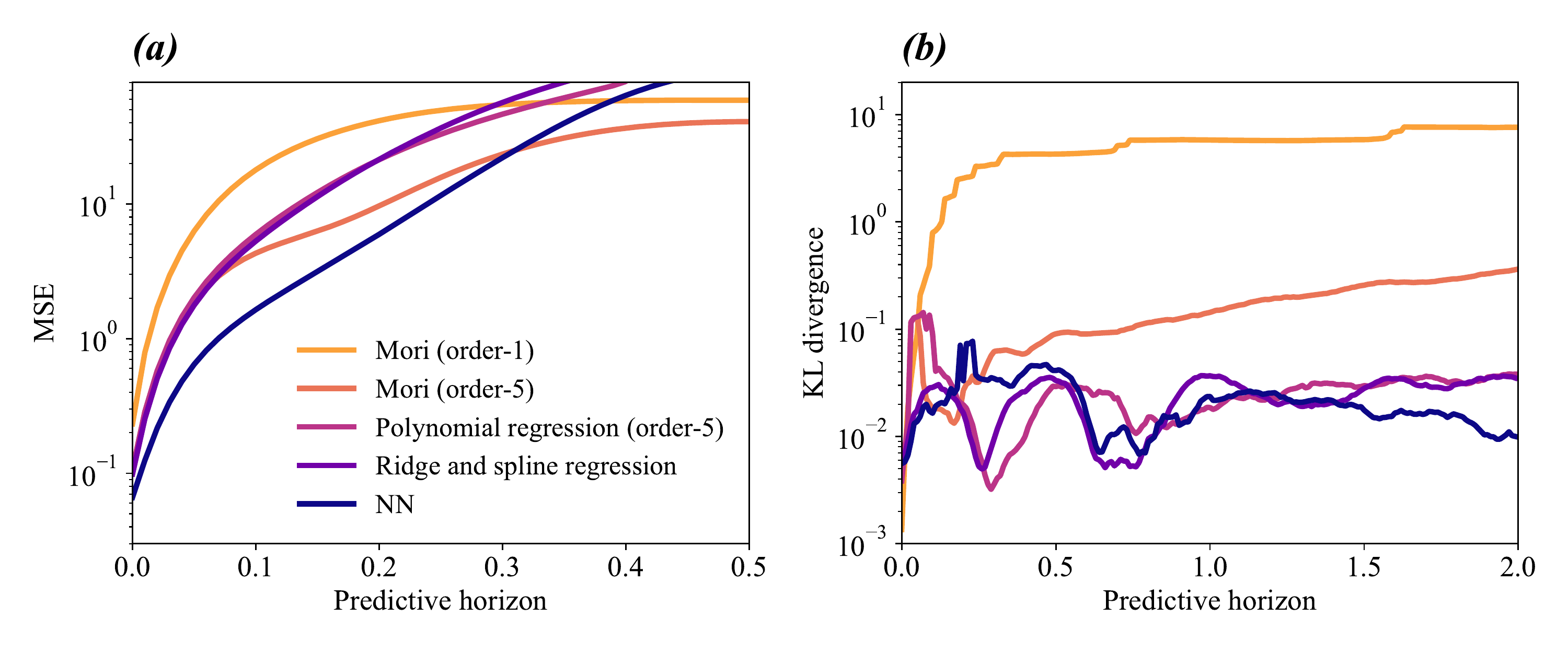}
    \caption{Error of the reduced-order models on the Lorenz (1963) system. (a) The mean squared error of the prediction as a function of the predictive horizon; (b) the Kullback--Leibler divergence of the empirical distribution from the model predicted distribution as a function of the predictive horizon.}
    \label{fig:Lorenz63Error} 
\end{figure}

We further quantified the prediction error by the following metrics. First, we computed the squared $L^2$-error of the prediction to the actual dynamics of $\varphi(t)$. The initial condition and the histories were uniformly sampled along a long trajectory of the actual dynamics. We collected 25,000 samples and compute the mean-squared error as a function of the predictive horizon in Fig.~\ref{fig:Lorenz63Error}(a). Within a short horizon, predictions by the nonlinear projections are better than linear projections. The errors of the nonlinear projections eventually grew larger at a finite predictive horizon ($t\approx 0.3$). We remark that, this should not be considered as a drawback but an advantage: beyond this timescale, linear projection predicted the mean of the chaotic dynamics, but nonlinear projection predicted oscillatory solutions. In terms of the mean squared error, predicting only the mean would have a smaller error. As such, we also computed the Kullback--Leibler divergence from the predictive distribution (approximated by the histogram of the samples from of the reduced-order models) to the ground truth distribution (approximated by the histogram of the samples from the actual dynamics), as shown in Fig.~\ref{fig:Lorenz63Error}(b). With this metric, we observed that the reduced-order models based on the nonlinear projection operators captured the time-dependent distribution of $\varphi$ much better than the reduced-order models based on linear projection operators. 

\begin{figure}[!t]
    \centering
    \includegraphics[width=\textwidth]{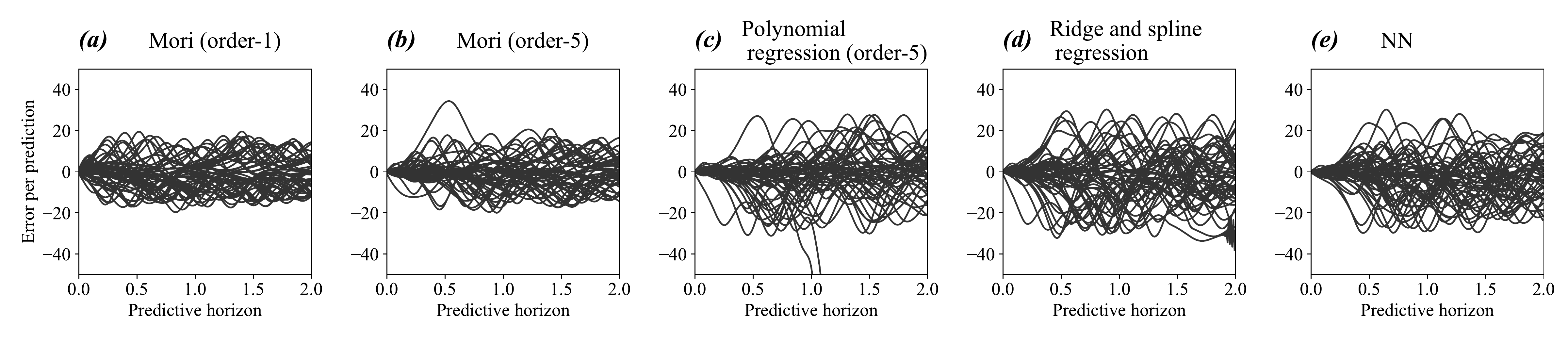}
    \caption{The Lorenz (1963) system. (a-e) The deviation of 50 predicted trajectories using the learned models.}
    \label{fig:Lorenz63TrajError} 
\end{figure}

We conclude this section by discussing the potential pitfalls of regression-based nonlinear projection operator. Mainly, the stability of the learned models is not always guaranteed. We observed that the learned reduced model, provided some initial histories, can blow up in finite time. In Fig.~\ref{fig:Lorenz63TrajError}, we plotted the deviation of the prediction from the actual dynamics. For each of the reduced-order models, we plotted 50 deviations, whose initial condition and histories sampled on a long ground-truth trajectory. We observed that the prediction the plain fifth-order polynomial regression could explode in finite time. The spline regression, which predicts constants outside the boundary knots seemed to alleviate the problem. However, when the prediction goes beyond the data distribution, the spline regression can also be ``trapped'' locally and induce high-frequency oscillations. Interestingly, we did not observe these pathologies in the neural network.
 
\subsection{Kuramoto--Sivashinsky Equation}\label{sec:KSmodel}
Our third example is the Kuramto--Sivashinsky (KS) equation, which was developed in the 1970s for modeling instability of flame propagation \cite{kuramotoDiffusioninducedChaosReaction1978,sivashinskyFlamePropagationConditions1980,sivashinskyNonlinearAnalysisHydrodynamic1977}. Specifically, in this study, we considered the one-dimensional KS (partial differential) equation
\begin{subequations}
\begin{align}
    0 = {}& \frac{\partial}{\partial t} u\l(x,t\r) + \frac{\partial^2}{\partial x^2}u\l(x,t\r) + \frac{\partial^4}{\partial x^4} u\l(x,t\r) + u\l(x,t\r) \frac{\partial}{\partial_x} u\l(x,t\r),
\end{align}
\end{subequations}
for all time $t\ge 0$ on a bounded domain $x\in\l[0,L\r]$, given initial data $u\l(x,0\r) = u_0\l(x\r)$. A periodic boundary condition was imposed, $u\l(t,0\r)=u\l(t,L\r)$ $\forall t$. We chose $L=16 \pi$ and uniformly discretized the space into $N=128$ points. We used the pseudo-spectral forth-order exponential time-derivative Runge--Kutta (ETDRK4) algorithm \cite{kassamFourthOrderTimeSteppingStiff2005} with 2/3-dealiasing \cite{linDatadrivenModelReduction2021} to solve the spatially discretized system. The time step for the explicit ETDRK4 algorithm was set as $10^{-3}$, and we collect the snapshots every $10^3$ steps (so $\Delta=1$). For the training data, we followed \cite{kassamFourthOrderTimeSteppingStiff2005} and set $u_0^{\text{tr}}(x)=\cos(2 \pi x/L)\l[1+\sin(2\pi x/L)\r]$, and for the test data, $u_0^{\text{te}}(x)=\sin(2 \pi x/L)\l[1+\cos(2\pi x/L)\r]$. Both the training and the validation sets were simulated to $t=t_c+10^5$, and we discarded the transient trajectories before $t_c=5\times 10^2$. For a reference, the Lyapunov time of the system is $\approx 12.56$ \cite{edson_bunder_mattner_roberts_2019}.  We still chose mean-square error as the cost function for the regression models. 

We considered a four-fold reduced-order model by only observing the variable $u$ every four discretized spatial grids. That is, $u\l(t,x_{4i}\r)_{i}$, where $i=0,1,\ldots 31$ is the grid index. We augmented the reduced-order data set by two operations. (1.)~\emph{Shifting}: we used samples collected on different sub-grids: $u\l(t,x_{4i+j}\r)$, $j=0,1,2,3$ and (2.)~\emph{Reordering}: we imposed the $C_1$ symmetry of the system, i.e., $u\l(t,x_{\text{mod}\l(4\l(i+k\r)+j, 32\r)}\r)$, $k=0,1,\ldots 31$. 

We also adopted and integrated the delay-embedding technique into some of the regression-based models. With the delay embedding, the input of the regression analysis is augmented to include a finite number, $E \in \mathbb{Z}_{\ge 0}$, of the current and past snapshots \cite{arbabiErgodicTheoryDynamic2017,luschDeepLearningUniversal2018}. As remarked in the discussion in our previous study \cite{lin2021DataDrivenLearningMori}, the history dependence in the delay-embedding technique is not the memory effect in the Mori--Zwanzig formalism. We will provide a more detailed discussion in the Discussion section. For the rest of this article, we  refer to the former as ``delay-embedding'' and the latter as ``(MZ) memory effect'' ``(MZ) memory correction'', or simply ``memory'' in the figures. Our major aim of the analysis is to compare the performance of the models with delay-embedding and without MZ memory correction to those models with MZ memory correction but without delay-embedding.

We considered the following regression-based models and compared their performance. The first model is the linear Mori projection operator. Our preliminary analysis showed that a direct application of either data-driven Koopman \cite{schmid2010DynamicModeDecomposition} or Mori's \cite{lin2021DataDrivenLearningMori} on the snapshot observations $\l\{u\l(\cdot, x_{4i}\r)\r\}_{i=0}^{31}$ performed poorly [data not reported]. Thus, we adopted the delay-embedding technique \cite{arbabiErgodicTheoryDynamic2017,luschDeepLearningUniversal2018} to augment the basis functions for learning the MZ kernels based on Mori's projection operator (i.e., linear projection with memory). The model will be referred to as ``Mori+DEm'' model. Specifically, we use the past $E$ snapshots to augment the input vector (which lives in $\mathbb{R}^{32E}$). For Mori's projection operator, we chose a long embedding $E=10$. We remark that the lowest-order result of the Mori's projection operator (i.e., without the higher-order memory corrections, $\boldsymbol{\Omega}^{(\ell)}$, $\ell\ge 1$, in Eq.~\eqref{eq:GLE})) is functionally identical to the Hankel-DMD \cite{arbabiErgodicTheoryDynamic2017}, which is a method solely based on delay embedding. The second regression-based model was chosen to be a fully-connected neural network (FCNN). The architecture of the NN was two fully-connected layers, each of which contains 32 artificial neurons and a third linear layer with 32 outputs. The third regression-based projection operator was a convolutional neural network (CNN). CNNs are translationally invariant; such inductive bias is considered as beneficial for systems like KS equation, which are invariant under translations. The architecture of the CNN was two one-dimensional convolutional layers with five channels and kernel size of 11. The circular padding was implemented to impose the periodic boundary condition. As the fourth projection operator, we adopted the delay-embedding technique and used $E=4$ past snapshots as the input for a CNN model, which will be referred as ``CNN+DEm'' model. In practice, the past $E$ snapshots were stacked in the channel dimension of the input tensor of CNN. The architectures of the delay-embedded CNN model is identical to the one without the delay-embedding. This model will be referred to as the ``CNN+DEm'' model.

We used the snapshots collected along the long trajectory with the initial condition $u^\text{tr}_0$ to learn the MZ kernels $\mathbf{\Omega}^{\l(\ell\r)}$ for each regression-based projection operators. For linear projection operators and for FCNN, we applied both the shifting and reordering data augmentation to impose the symmetries. For CNN's, we applied only the shifting data augmentation. We sampled from the generated test trajectory to perform error analyses of the learned regression-based MZ kernels.

\begin{figure}[!t]
    \centering
    \includegraphics[width=\textwidth]{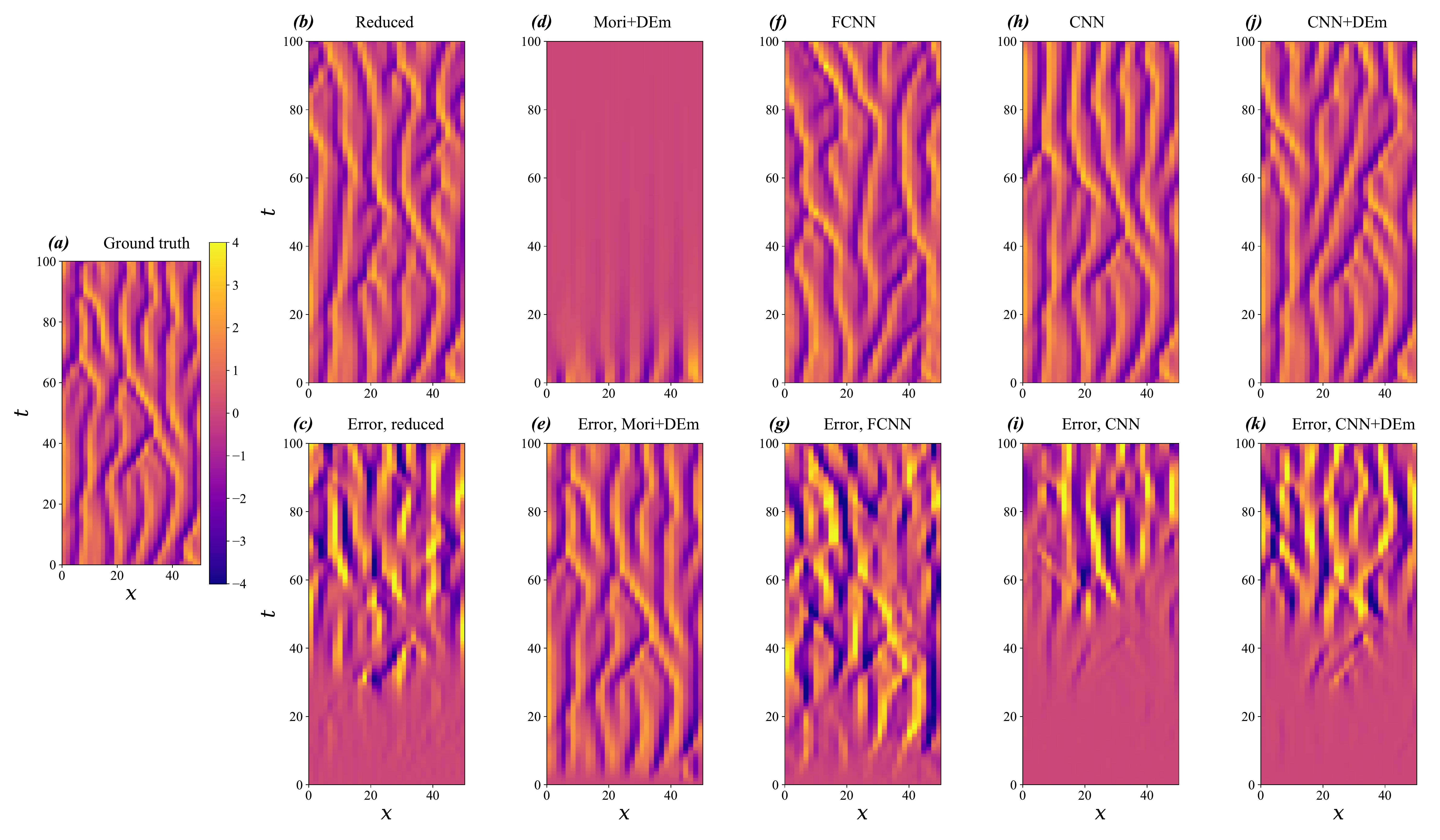}
    \caption{Regression-based Mori-Zwanzig learning of Kuramoto--Sivashinsky equation. (a) ``Ground-truth'' reduced-order observation. (b, d, f, h, j) Predicted trajectories using reduced-order simulation and four regression based Mori--Zwanzig models and (c, e, g, i, k) error of the predicted trajectories. The $x$- and $y$-axis of panels (b-k) are identical to those in the panel (a).}
    \label{fig:KSTrajs} 
\end{figure}

As a baseline reference for comparison, we performed reduced-order simulations. In these simulations, the reduced-order snapshots ($\l\{u_{t, x_{4i}}\r\}_{i=0}^{31}$) were treated as the full state of the KS-equation discretized on 32 spatially discretized grids (v.~128 of the ground-truth data-generation process). Each of the snapshots was evolved forward in time directly by the ETDRK4 integrator. 

\begin{figure}[!t]
    \centering
    \includegraphics[width=\textwidth]{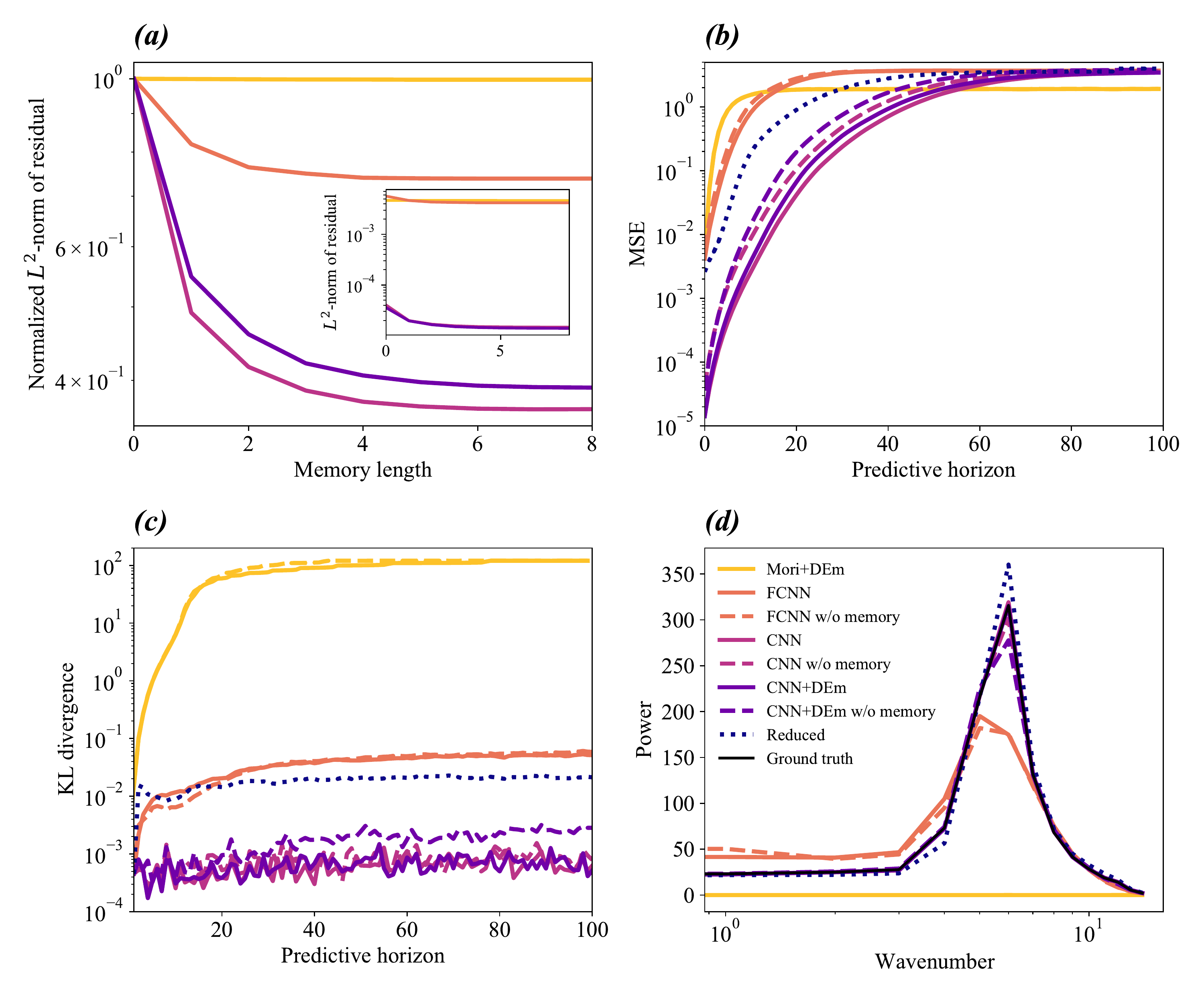}
    \caption{Mori--Zwanzig learning of Kuramoto--Sivashinsky system with regression-based projection operator. (a) The mean squared error of the learned models as a function of memory length. (b) The mean squared prediction error over different prediction horizon. The Lyapunov time for chosen domain size $L=16\pi$ is approximately 12.56. (c) The Kullback--Leibler divergence of the marginal ($u\l(t, x_{4i}\r)$ $\forall i=0\ldots 31$) ground-truth distribution from the predicted distribution as a function of prediction horizon. (d) The power spectrum of different models at a prediction horizon 100.}
    \label{fig:KSError} 
\end{figure}

In Figure \ref{fig:KSTrajs}, we visualized the predictive trajectories and errors from the learned models, each with $\ell \le 10$ MZ memory kernels. We chose $H=10$ and used the truncated GLE \eqref{eq:GLEtruncated} for making the predictions, again with the assumption $\mathbf{W}_n=0$. Figures \ref{fig:KSTrajs} (b, c) show the spatio-temporal evolution of the reduced-order simulation. Qualitatively speaking, we observed that the reduced trajectory deviated from the ``ground-truth'' at a prediction horizon $\approx 20$, but it is capable of reproducing dynamical features that are visually indistinguishable. Mori's projected model with $E=10$ delay-embedding (Figs.~\ref{fig:KSTrajs} (d, e)) quickly converged to the mean of the dynamics, as we expected. The FCNN model deviated from the ``ground-truth'' trajectory sooner than the reduced simulation did, but it also reproduced similar dynamical features; see Figs.~\ref{fig:KSTrajs} (f, g). Predictions from both CNN models deviated later than the reduced simulation did (Figs.~\ref{fig:KSTrajs} (h-k)). The $E=4$ delayed-embedded CNN (Figs.~\ref{fig:KSTrajs} (j, k)) did not significantly improve the plain CNN (Figs.~\ref{fig:KSTrajs} (h, i)).

Quantitative error analysis of the regression models were performed with 15,000 trajectories uniformly sampled on the generated test trajectory. Figure \ref{fig:KSError}(a) shows the normalized residual of the one-step prediction, defined as the mean squared error of the $\ell$th-order MZ model, normalized by that of the Markovian model (i.e., only $\boldsymbol{\Omega}^{\l(0\r)}$). The absolute error values are presented in the inset of Fig.~\ref{fig:KSError}(a). For all the regression models, including the MZ memory reduced the prediction error. We observed that the MZ memory effects can be captured with less than 10 memory kernels, justifying our choice of $H=10$ for the predictive model. The plain CNN model without delay-embedding showed the highest percentage improvement by including the MZ memory. The delay-embedded CNN model had a slightly smaller percentage improvement in comparison to that of the plain CNN model, despite it has smaller absolute magnitude of mean squared error. The smaller absolute error of the $E=4$ delay-embedded model is expected, as it is a more expressive model due to its larger input space, allowing coupling the delayed observables. Surprisingly, the plain CNN with $H=4$ MZ contribution has a smaller absolute error than the $E=4$ delay-embedded CNN (inset of Fig.~\ref{fig:KSError}(a)); see the Discussion Section for a detailed discussion. The FCNN model exhibited a much larger magnitude of error, and a smaller improvement from MZ memory contributions. Finally, as we expected, the linear Mori projection operator with delay-embedding had a very poor performance in one-step prediction.

\begin{figure}
    \centering
    \includegraphics[width=\textwidth]{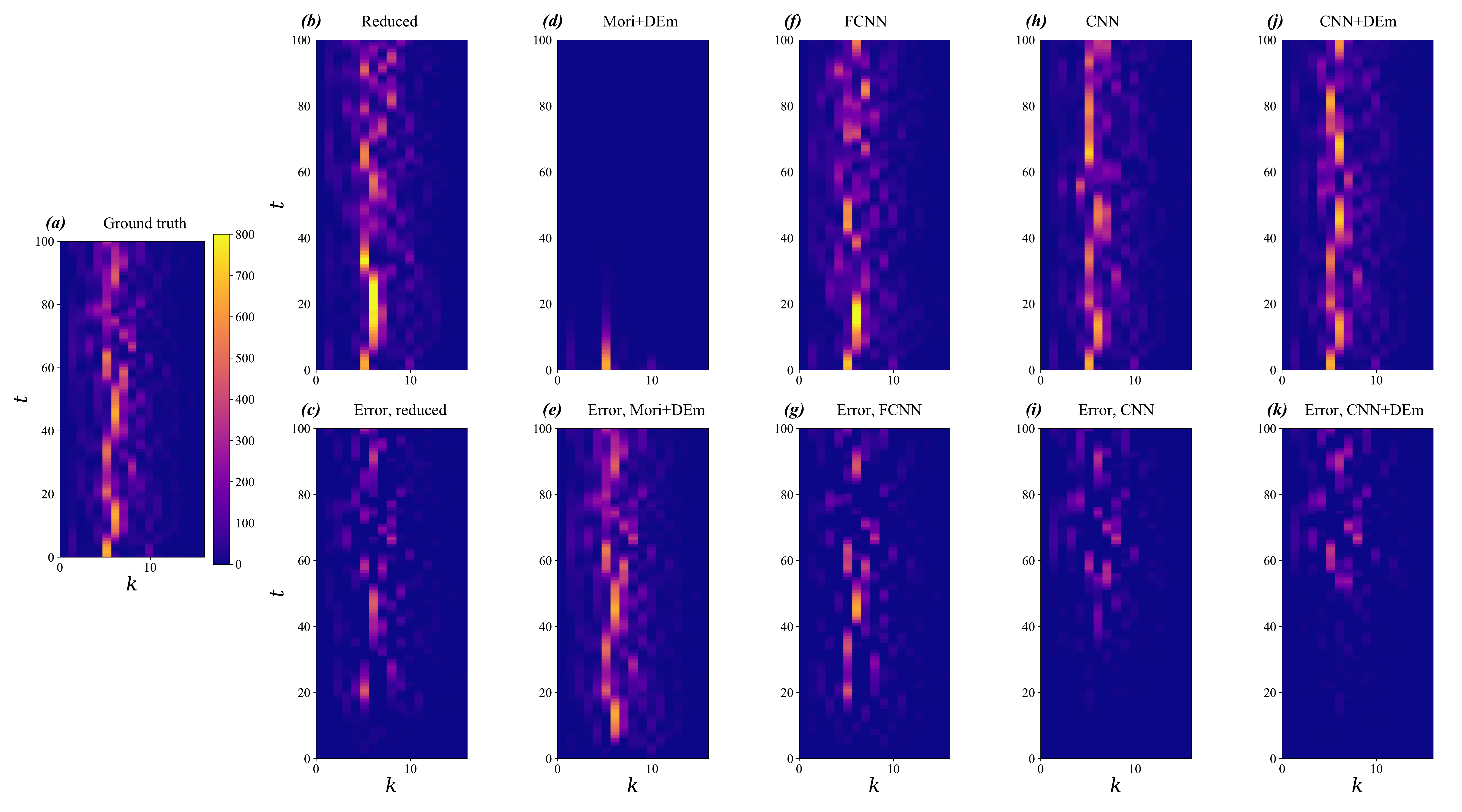}
    \caption{Projection-based Mori--Zwanzig learning of Kuramoto--Sivashinsky system. (a) ''Ground-truth`` spectrum of the coarse-grained trajectory. (b, d, f, h, j) Predicted spectrum using reduced-order simulation and four regression based Mori--Zwanzig models and (c, e, g, i, k) error of the predicted spectrum.}
    \label{fig:KSSpectrum} 
\end{figure}

We next focused on predictions at a further horizon. We set the furthest prediction horizon at $t=100$, which is approximately 8 Lyapunov time. In Figure \ref{fig:KSError} (b), we compared the mean squared error as a function of the predictive horizon. It was observed that including the MZ memory improved the multi-step predictions. Similar to the Lorenz ’63 system, at short prediction horizon ($\lesssim 20$), models with nonlinear projection operator performed much better than the linear Mori's model. Among the nonlinear models, the FCNN model had worse predictions than the reduced simulation, and both the CNN models out-performed the reduced simulation. Comparing both CNN models, the delay-embedded CNN models performed better than the plain CNN at one-step prediction, as established in Fig.~\ref{fig:KSError} (a). However, the CNN model overtook the time-embedding CNN model and became the best model---in the metric of the mean squared error---among all the considered models before the prediction horizon exceeded $55$.  As illustrated with the Lorenz (1963) model, the mean squared error is a less meaningful metric for the long-horizon predictions. As a second metric to quantify model performances, we used the Kullback--Leibler (KL) divergence from the regression-model predicted marginal distribution ($u$ on a single grid) to that of the ground-truth (test trajectory) distribution. Because of the translational symmetry, we accumulated all the measured $u\l(t, x_{4i}\r), i=0\ldots 31$, at different horizons for computing the empirical one-dimensional distributions of $u$. Figure \ref{fig:KSError} (c) shows the computed Kullback--Leibler divergence from the predicted distributions at different prediction horizon. Finally, we compared the power spectra of based  at the prediction horizon 100; Figure \ref{fig:KSSpectrum} also shows the spectrograms of the models with a single sampled initial condition. With these metrics, we concluded that the CNN models, with the MZ memory, performed much better than all other models and the reduced simulations over a long prediction horizon. Again, the presented evidence illustrated the advantage of nonlinear projection at predicting statistical properties of the dynamical system at a longer prediction horizon, despite the fact that sample-wise, the predicted trajectories deviate from the ground-truth trajectories.

\begin{figure}[!t]
    \centering
    \includegraphics[width=\textwidth]{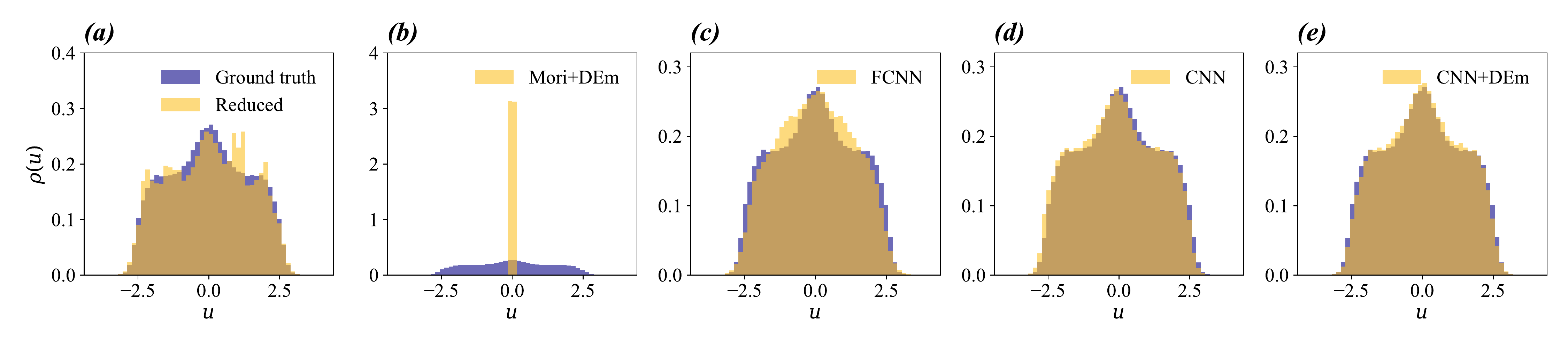}
    \caption{Projection-based Mori-Zwanzig learning of Kuramoto-Sivashinsky system. (a-e) PDFs of long-time predicted trajectories(horizon = 10,000) for five regression-based MZ models.}
    \label{fig:KSPDF} 
\end{figure}

The long-time marginal distribution of the regression-based models and the reduced simulation is shown in Figure \ref{fig:KSPDF}. The empirical distributions were computed by accumulating measured $u$ on each grid along long predicted trajectories (15,000 steps). The initial conditions of all the models, as well as the reduced simulation, were all set as $u^\text{te}_0$. Mori’s projection with delay-embedding showed almost $\delta$-like distribution at the mean as expected, while the nonlinear models were capable of reproducing the stationary marginal distribution similar to the ``ground-truth''. Among the nonlinear models, the distributions from both CNN models almost perfectly aligned with that of the test trajectory, while the FCNN model missed some finer details around the mean.

\section{Discussion and Conclusion}\label{sec:discussion}

The central proposition of this article is to adopt the regression analysis to define the projection operator in the Mori--Zwanzig formalism to enable data-driven learning of the operators. Given a set of sampled pairs of dependent and independent variables, a regression analysis identifies an optimal parametric function of the independent variables, such that the residuals are minimized in some norm, prescribed by the noise model. In the context of dynamical systems, the dependent variables are the future observations, and the independent variables are the present ones. When the system is not fully resolved, the future observations can not always be uniquely determined by the present ones due to missing information of the unresolved degrees of freedom. Consequently, future observations can be a one-to-many mapping of the present ones in under-resolved systems. A regression projects such one-to-many mapping to a one-to-one or many-to-one function, constrained by the data. It is a projection operator because an additional regression analysis on the samples drawn from the best-fit model results in the best-fit model, satisfying the definition of a projection operator $\hat{\mathcal{P}}$: $\hat{\mathcal{P}}^2=\hat{\mathcal{P}}$.

To identify regression analysis as a projection operator is not a novel idea; it is an established concept in statistics. For example, see \cite{deisenroth_faisal_ong_2020} for linear regression as an orthogonal projection, and \cite{nosedal2012reproducing} for regularized regression in terms of the Reproducing Kernel Hilbert Space. The notion of the decomposing the samples of the dependent variables into a regressional and an orthogonal component \cite{Afriat1969} even went back to Kolmogorov in the 1940s \cite{Kolmogorov46Regression}. The major conceptional contribution of this manuscript is the proposition \textit{to adopt a regression analysis as the projection operator in the Mori--Zwanzig formalism}. This is a novel concept which, to our best knowledge, does not exist in current literature. Then, we used the Generalized Fluctuation-Dissipation (GFD) relationship, a key self-consistent condition in the Mori--Zwanzig formalism, to recursively extract the memory kernels. Our novel approach enforces the GFD, which is neither enforced nor checked in existing methods \cite{weinanE,GILANI2021132829} motivated by the Mori--Zwanzig formalism. In doing so, we are able to directly extract the MZ memory operators computationally, without modeling them with additional assumptions. 

We showed that the special choice of adopting linear regression models results in our previously proposed data-driven learning algorithm for the MZ formalism with Mori's projection operator \cite{lin2021DataDrivenLearningMori}. Thus, the proposed regression-based projection operator is a generalization of our previous proposition.
Our proposed procedure to extract the MZ memory kernels is readily applicable to a large set of data-driven learning methods for dynamical systems in the literature because they can be cast into the form of Eqs.~\eqref{eq:generalCost}-\eqref{eq:regression-basedProjection}. For example, approximate Koopman methods with neural-network-based learning of the functional bases with ridge and sparsity-induced regularizations \cite{liExtendedDynamicMode2017,Yeung2017LearningDN,luschDeepLearningUniversal2018,wehmeyerTimelaggedAutoencodersDeep2018}, and Sparse Identification of Nonlinear Dynamics \cite{Brunton16SINDy,Kaheman20SINDyPI}. We remark that these data-driven learning methods often implicitly assume that the system is fully resolved without missing information. Should it be the case and should the regression model be expressive enough, our procedure would lead to a perfect Markov model $\boldsymbol{\Omega}^{(0)}$ that leaves zero residuals. Our interest, as well as the setup of the Mori--Zwanzig formalism, is on partially resolved dynamical systems. In general, when a system is not fully resolved and when there exists interaction between the resolved and unresolved dynamics, MZ memory kernels and orthogonal dynamics naturally emerge. 


Regression-based MZ methods can fill in the gap between the existing Mori's and Zwanzig's projection operators, which can be respectively considered as the simplest and the most optimal projection operator \cite{chorin2002OptimalPredictionMemory}. Our proposition significantly broadens the applicability of the Mori--Zwanzig formalism, because the complexity of the regression model can be gradationally adjusted between these two extremes. Ranked by the complexity of the regression analysis, we examined and presented the linear projection with \cite{lin2021DataDrivenLearningMori} and without memory \cite{schmid2010DynamicModeDecomposition,williams2015DataDrivenApproximation}, and various regression techniques such as polynomial regression, spline regression with ridge regularization, and the more modern and expressive neural network architectures such as FCNN and CNN for performing nonlinear projection, with or without memory. With a finite snapshot data set, it may not be computational feasible to learn the operators with the most optimal Zwanzig's projection operator, because it requires estimating the conditional expectation functions (e.g, $f(X):=\mathbb{E} \l[Y \vert X\r]$). A finite set of data may not be sufficient to cover all possible condition $X$, let alone the large amount of samples of $Y$ needed for estimating the expectation value. Instead, regression analysis identifies the best fit among a family of functions to approximate the conditional expectation function. When the family of functions is expressive enough to include the conditional expectation function, our regression-based MZ would converge to the Zwanzig's projection operator in the infinite-data limit. As neural networks are universal function approximators \cite{murphy2012machine}, we hypothesize that the MZ formalism with the optimal Zwanzig projection operator can be reasonably realized by a neural-network-based regression. Indeed, in our numerical experiments, we observed the superiority of those regression models based on neural networks.

Mathematically, the presented method learns the \emph{operators} $\boldsymbol{\Omega}^{(\ell)}$, which maps square integrable functions to square integrable functions, in Eq.~\eqref{eq:GLEo}. Computationally, this is made possible by using the regressional functions to model the \emph{functions} $\boldsymbol{\Omega}^{(\ell)}$ in Eq.~\eqref{eq:GLE}, given snapshots of the observations $\mathbf{g}\l(\boldsymbol{\phi}(t)\r)$, where the state is sampled from the induced distribution $\boldsymbol{\phi}\l(t\r)=\boldsymbol{\phi}\l(t; \boldsymbol{\phi}_0\r)$, $\boldsymbol{\phi}_0 \sim \mu$. It should be noted that (machine) learning of functional operators is not a new idea. For example, Li et al.~\cite{luLearningNonlinearOperators2021} recently laid out a general ML architecture (DeepONet) for learning nonlinear functional operators, relying on a universal approximation theorem proved in the 1990s \cite{tianpingchenUniversalApproximationNonlinear1995}. Another example is the approximate Koopman learning methods \cite{schmid2010DynamicModeDecomposition,Schmid2011,williams2015DataDrivenApproximation}, with which one aims to learn the (approximate) Koopman operator, linearly mapping functions to functions. With these methods, samples of function evaluations must be provided for supervised learning of the operators. Where the functions are evaluated---where the ``sensors'' are in the language of \cite{luLearningNonlinearOperators2021}---requires human inputs. With regard to our aim of learning dynamical systems, analogously, we also need to specify the initial distribution $\mu$. What is different from the general DeepONet \cite{luLearningNonlinearOperators2021} is that we can rely on the dynamics to naturally propagate the specified initial distribution to a future time $\boldsymbol{\phi}\l(t;\boldsymbol{\phi}_0\r)$, with which the samples of the observations can be generated for supervised learning---in our case, a general regression problem. 

We have been framing the above discussions in terms of data-driven learning methods. In fact, our original desire for developing data-driven MZ methods with nonlinear projection operators emerged from \emph{modeling} reduced-order dynamical systems. Devising closure schemes is not only necessary but also critical in these modeling endeavors, especially for multiscale models. Based on our knowledge, there are very few closure schemes which resemble the linear projection operators (i.e., approximate Koopman \cite{schmid2010DynamicModeDecomposition,williams2015DataDrivenApproximation} and Mori \cite{lin2021DataDrivenLearningMori}), with which one predicts the evolution of a set of observables by linear superposition of the same set of observables. More often, one adopts the nonlinear projection scheme defined in Sec.~\ref{sec:linearVnonlinear}, and iteratively apply the learned nonlinear map (i.e.~Eq.~\eqref{eq:mStepPredictionNonlinear}) to make multi-step predictions. By formally defining a regression (or more generally, data-driven parametrization) as the projection operator in the Mori--Zwanzig formalism, Algorithm \ref{alg:1} provides a principled way to extract the Mori--Zwanzig memory for a wide class of data-driven modeling and systems identification problems. Unless the residuals are all zero, the extracted memory contributions must be considered to account for the error due to the unresolved dynamics. Including the MZ memory in prediction reduces the error in prediction, as we observed in multiple examples provided in Sec.~\ref{sec:numerics}. 


In this article, we are primarily interested in the quality of multi-step predictions of the learned models. We showed that despite linear regression models out-performed nonlinear ones at the far horizon, it is only because the linear models predicted the mean of the dynamics. Because linear models fully ignore higher-order statistical properties of the dynamical systems, we are not confident that such ``coarse-grained dynamical systems'' are the proper way forward for making predictions. Although such a problem for linear projection operators does not exist when a set of observables linearly spanned an invariant Koopman space, it is extremely difficult to identify the invariant Koopman space in practice. It is our perspective that developing methods based on nonlinear projections (see Sec.~\ref{sec:linearVnonlinear}) is the way forward for predictive coarse-grained models. We do remark that methods based on linear projection operators are ideal for estimating the Koopman eigenvalues and eigenfunctions, and the globally convex learning problem is easier, provided a set of linearly independent observables. 
We should also remark that nonlinear regression is not always advantageous in making predictions. As we witnessed in the Lorenz (1963) model, despite a learned nonlinear regression model could improve predictive accuracy within a finite horizon, sometimes these nonlinear closure methods loose long-term stability (Fig.~\ref{fig:Lorenz63Error}). Such a pathology is model- and horizon-dependent, and interestingly, not observed when we used neural-network-based regression. Enforcing numerical stability in the process of learning  merits further investigation and may be a fruitful direction in the future. 
 
With the Kuramoto--Sivashinsky model, we also compared the delay-embedding and the memory effect of the Mori--Zwanzig formalism. Our motivation was to provide clarification on a pervasive confusion between these two memory-dependent theories. Such a confusion emerged from a false presumption that all memory-dependent theories are equivalent. Here, we elaborate the differences between these two approaches.

Because the memory effect in Mori--Zwanzig formalism is due to the interaction between the resolved and unresolved dynamics, the memory kernels and the orthogonal dynamics must satisfy the stringent Generalized Fluctuation Dissipation relation (Eq.~\eqref{eq:GFD}), which depends on the choice of the resolved observables and the projection operator. Guided by the GFD \eqref{eq:GFD}, learning MZ operators can be decomposed into $H$ smaller regression (optimization) problem in $M$ dimensions. Notably, the MZ memory contribution is always a linear sum of memory kernels evaluated at each past snapshot (Eq.~\eqref{eq:GLE}). The formalism excludes coupling between the observations made at different times, for example, $\mathbf{g}_{n_1} \times \mathbf{g}_{n_2}$, with $n_1\ne n_2$. This is because the projection operator $\mathcal{P}$ is applied at each time by construction. The functional space that $\mathcal{P}$ is projected into is those regressional functions $f\l(\cdot;\theta\r)$, which always only depend on $M$ reduced-order observables evaluated at a specific time. 

On the contrary, for delay-embedding, the domain of the functions is augmented to include $E$ past snapshots of $M$ reduced-order observables. The learning is a global regression (optimization) problem in an $M\times E$-dimensional functional space. The motivation of delay-embedding is that the past history may contain useful information for inferring the unresolved dynamics. Because such an augmented functional space ($M\times E$) is larger than the fixed ($M$) operational space in the Mori--Zwanzig formalism, delay-embedding can be considered as a more expressive regression model. Let us denote the past $E$ snapshots of the observables $\mathbf{g}_{0,-1\ldots -E+1}$ by a flattened vector $\mathbf{g}_{0:-E+1}\l(\boldsymbol{\phi}_{-T}\r) \in \mathbb{R}^{ME}$. Formally, delay embedding can be considered as a regression to identify the optimal function $\mathbf{f}_{\boldsymbol{\theta}^{\ast}}\l(\mathbf{g}_{0:-E+1} \l(\boldsymbol{\phi}_{-T}\r) \r)$ that best approximates $\mathbf{g}_{1}\l(\boldsymbol{\phi}_0\r)$.
As such, coupled terms could exist, should the family of regressional functions admits, within the delay-embedding paradigm. Hankel-DMD \cite{arbabiErgodicTheoryDynamic2017} is the special delay-embedding case that the learning is formulated as a linear regression problem, i.e., the regression functions take the form $\mathbf{f}\l(\mathbf{g}_{0:-E+1}\l(\boldsymbol{\phi}_0\r);\boldsymbol{\beta} \r):= \boldsymbol{\beta} \cdot \mathbf{g}_{t-E+1:t}\l(\boldsymbol{\phi}_0\r)$, where $\boldsymbol{\beta}$ is an $(ME)\times(ME)$ matrix (v.~in linear regression-based MZ, $\boldsymbol{\beta}^{(\ell)}$ is $M\times M$, $\ell=0\ldots H-1$ with the same memory length). Although its structure resembles Mori--Zwanzig's Generalized Langevin Equation, it is not necessary---and generally not the case---that the best-fit parameters of the Hankel-DMD satisfy the Generalized Fluctuation-Dissipation relation \ref{eq:GFD}. Notably, an extreme of the delay embedding technique is the Wiener projection, which has an infinitely long delay embedding. For stationary processes, Wiener projection leads to no Mori--Zwanzig memory \cite{linDatadrivenModelReduction2021}. A similar notion was also concluded in \cite{GILANI2021132829}, in which the authors argue that kernel methods could be used to learn the delay embedding with sufficient information that fully resolves the dynamics, resulting in zero Mori--Zwanzig memory contribution. Our proposition complements these approaches which attempt to eliminate the Mori--Zwanzig memory kernels: given any projection operator, our proposed procedure extracts the Mori-Zwanzig memory kernels.

\begin{figure}
    \centering
    \includegraphics[width=0.5\textwidth]{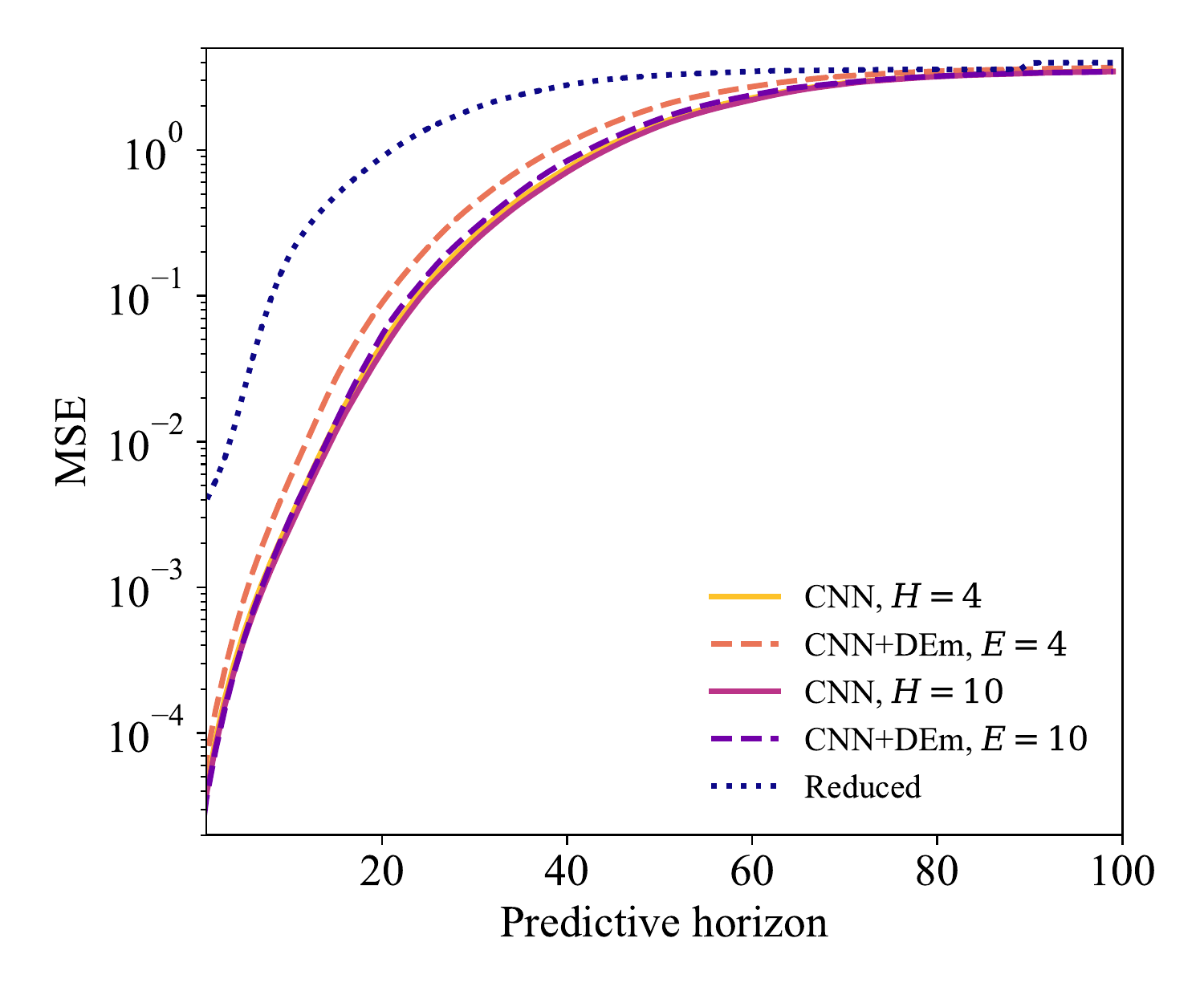}
    \caption{Mori--Zwanzig learning of Kuramoto--Sivashinsky system with various CNN-based regression model. The mean squared error of the prediction as a function of the prediction horizon. The CNN-based MZ models with different MZ memory lengths ($H=4$ and $H=10$) are compared with corresponding delay-embedded models with identical set of embedding lengths ($E=4$ and $E=10$).  The CNN architectures for this numerical experiment are specially designed such that: (1) the total number of parameters of both CNN-based MZ model and CNN-based delay-embedded model are the same and (2) the same information, i.e.~the length of the past history, is used for making prediction for both models (H=E).}
    \label{fig:KSSupL2KL} 
\end{figure}

Because the structure of the Generalized Langevin Equation is contained in the family of delay-embedded regressional functions, the error of the former should have been lower-bounded by that of the latter. However, with the Kuramoto--Sivashinsky equation, we showed that the Mori--Zwanzig formalism with $H = 4$ truncated memory contribution out-performed the more expressive delay-embedded CNN with $E=4$. We hypothesized that it was because the architecture of the delay-embedded CNN was not as expressive as the CNN with MZ memory. As the latter contained four separate CNNs, each of which learned an operator $\boldsymbol{\Omega}^{(\ell)}$, the delay-embedded CNN had a fixed architecture (the past history augmented the channel dimension of the input and not the complexity of the regression model). To examine whether a more expressive delay-embedded CNN could learn better, we deliberately developed delay-embedded CNNs whose number of parameters are identical to $H$ times the number of parameters CNNs for learning MZ operators. For $H=E=4$, we still observed that the MZ out-performed the delay-embedding model (1-step MSE of the CNN+MZ: $1.563\times 10^{-5}$; 1-step MSE of the CNN+DEm: $2.038\times 10^{-5}$). For $H=E=10$, the delay-embedding model performed better in predicting the first step than the MZ did (1-step MSE of the CNN+MZ: $1.478 \times 10^{-5}$; 1-step MSE of the CNN+DEm: $1.085\times 10^{-5}$). However, for multi-step predictions, MZ out-performed the delay-embedding after a few ($\mathcal{O}(10)$) steps, see Fig.~\ref{fig:KSSupL2KL}. These observations highlight the practical difficulty of learning the larger delay-embedded models: they require more data, potentially need a longer history, and for nonlinear problems, the training could be trapped in some local minimum. The precise reason why the delay-embedded models are worse in predicting multiple steps into the future merits a future, more carefully designed investigation. 

Furthermore, we demonstrated that Mori--Zwanzig formalism and delay embedding are not mutually exclusive. On the Kuramoto--Sivashinsky model, we illustrated a delay-embedded CNN combined with the MZ formalism. In this combined approach, the input of the Markov operator was a flattened observable $\mathbf{g}_{0:-E+1}$, and the input of the $\ell$th-order ($\ell \ge 1$) memory operator was a vector function of the flattened observable shifted by $\ell$ time step, $\mathbf{g}_{-\ell:-E-\ell+1}$. We showed in Fig.~\ref{fig:KSError}(a) and (b) that the inclusion of the MZ memory operators can further improve the performance of the delay-embedded CNN. 

We view regression-based Mori--Zwanzig formalism as a promising direction to enable data-driven learning for partially observed dynamical systems. We have identified a few directions which merit further investigation. Theoretically, we would like to generalize the formalism for stochastic systems, including generic stochastic processes, Hidden Markov Models, and random dynamical systems.
As an application, we are working on the application of the NN-based Mori--Zwanzig models on isotropic turbulence data \cite{tian2021DatadrivenLearningMori}. Another interesting idea we would like to explore is making formal connections to time-series analysis. As recently pointed out, Dynamic Mode Decomposition can be considered as a Vector Autoregression (VAR-1) model in time-series analysis \cite{bolltExplainingSurprisingSuccess2021}. It is clear that the Hankel-DMD can be mapped to a VAR-$p$ model in time-series analysis. It was shown that the Wiener projection can be mapped to the nonlinear autoregressive moving average model with exogenous inputs (NARMAX) in time-series analysis, after a rational approximation in the Laplace domain \cite{linDatadrivenModelReduction2021}. These existing evidence prompted an interesting research question: what is the time-series analysis method which corresponds to the MZ formalism, which imposes a very stringent GFD condition on the residuals and the memory? 
Next, we discussed the difference between the memory-embedding techniques such as Hankel-DMD and the MZ formalism. These two formalisms can be considered as two extremes of memory-dependent dynamics: The former considers an extremely large functional space which couples a whole segment of the past trajectory, and the latter restricts its functional space to the regressional functions at an instantaneous time. In modern machine learning for time series, Recurrent Neural Networks with Long Short-Term Memory (RNN-LSTM) \cite{HochreiterS97Long} and its variants are also equipped with memory kernels, learned by a neural network. The architecture has been claimed to learn Mori--Zwanzig memory kernels \cite{weinanE,HARLIM2021109922}. We argue that because the critical GFD is neither enforced nor checked, it is unclear that RNN-LSTM is learning MZ memory, or some other form of memory-dependent predictors. It will be an interesting direction to perform a careful analysis, dissecting the RNN-LSTM to identify what type of memory-dependent predictor RNN-LSTM is.
Finally, in this manuscript, we did not attempt to model the orthogonal dynamics, which is needed for making accurate predictions. Our test showed that the common choice of white (independent in time) Gaussian noise is unsatisfactory (data not shown). A certain color-noise must be needed. This is because the orthogonal dynamics, which encode unresolved dynamics, are also correlated in time. Despite that our numerical extraction of the orthogonal dynamics provides the data stream for learning, how one can use these extracted data to model orthogonal dynamics in a principled and system-agnostic way remains an open and challenging problem.

\section*{Acknowledgement}{YTL was primarily supported by the Laboratory Directed Research and Development (LDRD) Project ``Uncertainty Quantification for Robust Machine Learning'' (20210043DR). YT, DP and YTL were also supported by the LDRD project ``Accelerated Dynamics Across Computational and Physical Scales'' (20220063DR). DL was partially supported by the LDRD project ``MASS-APP: Multi-physics Adaptive Scalable Simulator for Applications in Plasma Physics'' (20220104DR).}

\bibliographystyle{vancouver}
\bibliography{sample}

\begin{thebibliography}{10}

\bibitem{mori1965transport}
Mori H.
\newblock Transport, Collective Motion, and {{Brownian}} Motion.
\newblock Progress of theoretical physics. 1965;33(3):423-55.

\bibitem{zwanzig1973nonlinear}
Zwanzig R.
\newblock Nonlinear Generalized {{Langevin}} Equations.
\newblock Journal of Statistical Physics. 1973;9(3):215-20.

\bibitem{chorin2002OptimalPredictionMemory}
Chorin AJ, Hald OH, Kupferman R.
\newblock Optimal Prediction with Memory.
\newblock Physica D: Nonlinear Phenomena. 2002 Jun;166(3-4):239-57.

\bibitem{parishDynamicSubgridScale2017}
Parish EJ, Duraisamy K.
\newblock A Dynamic Subgrid Scale Model for {{Large Eddy Simulations}} Based on
  the {{Mori}}\textendash{{Zwanzig}} Formalism.
\newblock Journal of Computational Physics. 2017 Nov;349:154-75.

\bibitem{parishNonMarkovianClosureModels2017}
Parish EJ, Duraisamy K.
\newblock Non-{{Markovian}} Closure Models for Large Eddy Simulations Using the
  {{Mori}}\textendash{{Zwanzig}} Formalism.
\newblock Phys Rev Fluids. 2017 Jan;2(1):14604.

\bibitem{stinis19}
Price J, Stinis P.
\newblock Renormalized Reduced Order Models with Memory for Long Time
  Prediction.
\newblock Multiscale Modeling \& Simulation. 2019;17(1):68-91.

\bibitem{stinis21}
Price J, Meuris B, Shapiro M, Stinis P.
\newblock Optimal renormalization of multiscale systems.
\newblock Proceedings of the National Academy of Sciences.
  2021;118(37):e2102266118.
\newblock Available from:
  \url{https://www.pnas.org/doi/abs/10.1073/pnas.2102266118}.

\bibitem{okamura2006validity}
Okamura M.
\newblock Validity of the essential assumption in a projection operator method.
\newblock Physical Review E. 2006;74(4):046210.

\bibitem{mori2007dynamic}
Mori H, Okamura M.
\newblock Dynamic structures of the time correlation functions of chaotic
  nonequilibrium fluctuations.
\newblock Physical Review E. 2007;76(6):061104.

\bibitem{meyer2020non}
Meyer H, Pelagejcev P, Schilling T.
\newblock Non-{M}arkovian out-of-equilibrium dynamics: A general numerical
  procedure to construct time-dependent memory kernels for coarse-grained
  observables.
\newblock EPL (Europhysics Letters). 2020;128(4):40001.

\bibitem{maeyama2020extracting}
Maeyama S, Watanabe TH.
\newblock Extracting and modeling the effects of small-scale fluctuations on
  large-scale fluctuations by {M}ori--{Z}wanzig projection operator method.
\newblock J Phys Soc Jpn. 2020;89(2):024401.

\bibitem{meyer2021numerical}
Meyer H, Wolf S, Stock G, Schilling T.
\newblock A numerical procedure to evaluate memory effects in non-equilibrium
  coarse-grained models.
\newblock Adv Theory Simul. 2021;4(4):2000197.

\bibitem{lin2021DataDrivenLearningMori}
Lin YT, Tian Y, Livescu D, Anghel M.
\newblock Data-{{Driven Learning}} for the {{Mori--Zwanzig Formalism}}: {{A
  Generalization}} of the {{Koopman Learning Framework}}.
\newblock SIAM Journal on Applied Dynamical Systems. 2021 Jan;20(4):2558-601.

\bibitem{tian2021DatadrivenLearningMori}
Tian Y, Lin YT, Anghel M, Livescu D.
\newblock Data-Driven Learning of {{Mori}}\textendash{{Zwanzig}} Operators for
  Isotropic Turbulence.
\newblock Physics of Fluids. 2021 Dec;33(12):125118.

\bibitem{rowley2009SpectralAnalysisNonlinear}
Rowley CW, Mezi{\'c} I, Bagheri S, Schlatter P, Henningson DS.
\newblock Spectral Analysis of Nonlinear Flows.
\newblock Journal of Fluid Mechanics. 2009;641:115-27.

\bibitem{schmid2010DynamicModeDecomposition}
Schmid PJ.
\newblock Dynamic Mode Decomposition of Numerical and Experimental Data.
\newblock Journal of Fluid Mechanics. 2010;656:5-28.

\bibitem{Schmid2011}
Schmid PJ, Li L, Juniper MP, Pust O.
\newblock Applications of the Dynamic Mode Decomposition.
\newblock Theoretical and Computational Fluid Dynamics. 2011 Jun;25(1):249-59.

\bibitem{williams2015DataDrivenApproximation}
Williams MO, Rowley CW, Kevrekidis IG.
\newblock A Data\textendash Driven Approximation of the {{Koopman}} Operator :
  Extending Dynamic Mode Decomposition.
\newblock Journal of Nonlinear Science. 2015;25(6):1307-46.

\bibitem{Durbin18}
Durbin P.
\newblock Some recent developments in turbulence closure modeling.
\newblock Annual Review of Fluid Mechanics. 2018;50:77-103.

\bibitem{Majda06}
Majda AJ, Wang X.
\newblock Nonlinear dynamics and statistical theories for basic geophysical
  flows.
\newblock {Cambridge University Press}; 2006.

\bibitem{Brunton16SINDy}
Brunton SL, Proctor JL, Kutz JN.
\newblock Discovering governing equations from data by sparse identification of
  nonlinear dynamical systems.
\newblock Proceedings of the National Academy of Sciences. 2016;113(15):3932-7.

\bibitem{liExtendedDynamicMode2017}
Li Q, Dietrich F, Bollt EM, Kevrekidis IG.
\newblock Extended Dynamic Mode Decomposition with Dictionary Learning: {{A}}
  Data-Driven Adaptive Spectral Decomposition of the {{Koopman}} Operator.
\newblock Chaos: An Interdisciplinary Journal of Nonlinear Science. 2017
  Oct;27(10):103111.

\bibitem{Yeung2017LearningDN}
Yeung E, Kundu S, Hodas NO.
\newblock Learning Deep Neural Network Representations for {{Koopman}}
  Operators of Nonlinear Dynamical Systems.
\newblock 2019 American Control Conference (ACC). 2017:4832-9.

\bibitem{luschDeepLearningUniversal2018}
Lusch B, Kutz JN, Brunton SL.
\newblock Deep Learning for Universal Linear Embeddings of Nonlinear Dynamics.
\newblock Nature Communications. 2018;9(1).

\bibitem{wehmeyerTimelaggedAutoencodersDeep2018}
Wehmeyer C, No{\'e} F.
\newblock Time-Lagged Autoencoders: {{Deep}} Learning of Slow Collective
  Variables for Molecular Kinetics.
\newblock The Journal of Chemical Physics. 2018;148(24):241703.

\bibitem{weinanE}
Model Reduction with Memory and the Machine Learning of Dynamical Systems.
\newblock Communications in Computational Physics. 2018;25(4):947-62.

\bibitem{HARLIM2021109922}
Harlim J, Jiang SW, Liang S, Yang H.
\newblock Machine learning for prediction with missing dynamics.
\newblock Journal of Computational Physics. 2021;428:109922.

\bibitem{GILANI2021132829}
Gilani F, Giannakis D, Harlim J.
\newblock Kernel-based prediction of non-Markovian time series.
\newblock Physica D: Nonlinear Phenomena. 2021;418:132829.

\bibitem{linDatadrivenModelReduction2021}
Lin KK, Lu F.
\newblock Data-Driven Model Reduction, {{Wiener}} Projections, and the
  {{Koopman-Mori-Zwanzig}} Formalism.
\newblock Journal of Computational Physics. 2021 Jan;424:109864.

\bibitem{chorinLu15}
Chorin AJ, Lu F.
\newblock Discrete approach to stochastic parametrization and dimension
  reduction in nonlinear dynamics.
\newblock Proceedings of the National Academy of Sciences. 2015;112(32):9804-9.

\bibitem{QIAN2020132401}
Qian E, Kramer B, Peherstorfer B, Willcox K.
\newblock Lift \& Learn: Physics-informed machine learning for large-scale
  nonlinear dynamical systems.
\newblock Physica D: Nonlinear Phenomena. 2020;406:132401.

\bibitem{gerlich1973VerallgemeinerteLiouvilleGleichung}
Gerlich G.
\newblock Die Verallgemeinerte {{Liouville}}-{{Gleichung}}.
\newblock Physica. 1973 Nov;69(2):458-66.

\bibitem{zwanzig2001nonequilibrium}
Zwanzig R.
\newblock Nonequilibrium Statistical Mechanics.
\newblock {Oxford ; New York}: {Oxford University Press}; 2001.

\bibitem{darveComputingGeneralizedLangevin2009}
Darve E, Solomon J, Kia A.
\newblock Computing Generalized {{Langevin}} Equations and Generalized
  {{Fokker-Planck}} Equations.
\newblock Proceedings of the National Academy of Sciences. 2009
  Jul;106(27):10884-9.

\bibitem{koopma32dynamical}
Koopman BO, v~Neumann J.
\newblock Dynamical Systems of Continuous Spectra.
\newblock Proceedings of the National Academy of Sciences. 1932;18(3):255-63.

\bibitem{arbabiErgodicTheoryDynamic2017}
Arbabi H, Mezi{\'c} I.
\newblock Ergodic Theory, Dynamic Mode Decomposition, and Computation of
  Spectral Properties of the {{Koopman}} Operator.
\newblock SIAM Journal on Applied Dynamical Systems. 2017;16(4):2096-126.

\bibitem{arbabi2018data}
Arbabi H, Korda M, Mezi{\'c} I.
\newblock A data-driven koopman model predictive control framework for
  nonlinear partial differential equations.
\newblock In: 2018 IEEE Conference on Decision and Control (CDC). IEEE; 2018.
  p. 6409-14.

\bibitem{yeung2019learning}
Yeung E, Kundu S, Hodas N.
\newblock Learning deep neural network representations for Koopman operators of
  nonlinear dynamical systems.
\newblock In: 2019 American Control Conference (ACC). IEEE; 2019. p. 4832-9.

\bibitem{Kaheman20SINDyPI}
Kaheman K, Kutz JN, Brunton SL.
\newblock SINDy-PI: a robust algorithm for parallel implicit sparse
  identification of nonlinear dynamics.
\newblock Proceedings of the Royal Society A: Mathematical, Physical and
  Engineering Sciences. 2020;476(2242):20200279.

\bibitem{fasel2022ensemble}
Fasel U, Kutz JN, Brunton BW, Brunton SL.
\newblock Ensemble-SINDy: Robust sparse model discovery in the low-data,
  high-noise limit, with active learning and control.
\newblock Proceedings of the Royal Society A. 2022;478(2260):20210904.

\bibitem{lorenzDeterministicNonperiodicFlow1963}
Lorenz EN.
\newblock Deterministic Nonperiodic Flow.
\newblock Journal of Atmospheric Sciences. 1963;20(2):130-41.

\bibitem{takens81detecting}
Takens F.
\newblock Detecting strange attractors in turbulence.
\newblock In: Rand D, Young LS, editors. Dynamical Systems and Turbulence,
  Warwick 1980. Berlin, Heidelberg: Springer Berlin Heidelberg; 1981. p.
  366-81.

\bibitem{kuramotoDiffusioninducedChaosReaction1978}
Kuramoto Y.
\newblock Diffusion-Induced Chaos in Reaction Systems.
\newblock Progress of Theoretical Physics Supplement. 1978 Feb;64:346-67.

\bibitem{sivashinskyFlamePropagationConditions1980}
Sivashinsky GI.
\newblock On Flame Propagation under Conditions of Stoichiometry.
\newblock SIAM Journal on Applied Mathematics. 1980;39(1):67-82.

\bibitem{sivashinskyNonlinearAnalysisHydrodynamic1977}
Sivashinsky GI.
\newblock Nonlinear Analysis of Hydrodynamic Instability in Laminar
  Flames\textemdash{{I}}. {{Derivation}} of Basic Equations.
\newblock Acta Astronautica. 1977;4(11):1177-206.

\bibitem{kassamFourthOrderTimeSteppingStiff2005}
Kassam AK, Trefethen LN.
\newblock Fourth-{{Order Time-Stepping}} for {{Stiff PDEs}}.
\newblock SIAM Journal on Scientific Computing. 2005 Jan;26(4):1214-33.

\bibitem{edson_bunder_mattner_roberts_2019}
Edson RA, Bunder JE, Mattner TW, Roberts AJ.
\newblock Lyapunov Exponents of the {{Kuramoto}}\textendash{{Sivashinsky PDE}}.
\newblock The ANZIAM Journal. 2019;61(3):270-85.

\bibitem{deisenroth_faisal_ong_2020}
Deisenroth MP, Faisal AA, Ong CS.
\newblock Mathematics for Machine Learning.
\newblock Cambridge University Press; 2020.

\bibitem{nosedal2012reproducing}
Nosedal-Sanchez A, Storlie CB, Lee TC, Christensen R.
\newblock Reproducing kernel Hilbert spaces for penalized regression: A
  tutorial.
\newblock The American Statistician. 2012;66(1):50-60.

\bibitem{Afriat1969}
Afriat SN.
\newblock In: Regression and Projection. Berlin, Heidelberg: Springer Berlin
  Heidelberg; 1969. p. 277-301.

\bibitem{Kolmogorov46Regression}
Kolmogorov AN.
\newblock On the proof of the method of least squares.
\newblock Uspekhi Mat~Nauk. 1946;1(1 (11)):57-70.

\bibitem{murphy2012machine}
Murphy KP.
\newblock Machine Learning: A Probabilistic Perspective.
\newblock {MIT press}; 2012.

\bibitem{luLearningNonlinearOperators2021}
Lu L, Jin P, Pang G, Zhang Z, Karniadakis GE.
\newblock Learning Nonlinear Operators via {{DeepONet}} Based on the Universal
  Approximation Theorem of Operators.
\newblock Nature Machine Intelligence. 2021 Mar;3(3):218-29.

\bibitem{tianpingchenUniversalApproximationNonlinear1995}
{Tianping Chen}, {Hong Chen}.
\newblock Universal Approximation to Nonlinear Operators by Neural Networks
  with Arbitrary Activation Functions and Its Application to Dynamical Systems.
\newblock IEEE Transactions on Neural Networks. 1995 Jul;6(4):911-7.

\bibitem{bolltExplainingSurprisingSuccess2021}
Bollt E.
\newblock On Explaining the Surprising Success of Reservoir Computing
  Forecaster of Chaos? {{The}} Universal Machine Learning Dynamical System with
  Contrast to {{VAR}} and {{DMD}}.
\newblock Chaos: An Interdisciplinary Journal of Nonlinear Science. 2021
  Jan;31(1):013108.

\bibitem{HochreiterS97Long}
Hochreiter S, Schmidhuber J.
\newblock Long Short-Term Memory.
\newblock Neural Comput. 1997;9(8):1735-80.
\newblock Available from: \url{https://doi.org/10.1162/neco.1997.9.8.1735}.

\end{thebibliography}

\end{document}